\documentclass[smallextended]{svjour3}
\newtheorem{Theorem}{Theorem}
\newtheorem{Lemma}[Theorem]{Lemma}

\newtheorem{Proposition}[Theorem]{Proposition}
\newtheorem{Definition}[Theorem]{Definition}
\newtheorem{Example}[Theorem]{Example}

\begin{document}
\title{Random planar graphs with bounds on the maximum and minimum degrees}
\author{Chris Dowden}
\institute{LIX, \'Ecole Polytechnique, 91128 Palaiseau Cedex, France \\
\email{dowden@lix.polytechnique.fr}}
\journalname{Graphs and Combinatorics}
\date{June 18, 2010}
\maketitle

\begin{abstract}

Let $P_{n,d,D}$ denote the graph taken uniformly at random from the set
of all labelled planar graphs on $\{ 1,2, \ldots, n \}$
with minimum degree at least $d(n)$ and maximum degree at most $D(n)$.
We use counting arguments to investigate the probability that $P_{n,d,D}$
will contain given components and subgraphs,
showing exactly when this is bounded away from $0$ and $1$ as $n \to \infty$.

\end{abstract}

\section{Introduction}

Random planar graphs have recently been the subject of much activity,
and many properties of the standard random planar graph $P_{n}$
(taken uniformly at random from the set of all labelled planar graphs on $\{1,2, \ldots, n\}$)
are now known.
For example,
in~\cite{mcd} it was shown that $P_{n}$ will asymptotically almost surely
(that is, with probability tending to $1$ as $n$ tends to infinity)
contain at least linearly many copies of any given planar graph.
By combining the counting methods of~\cite{mcd}
with some rather precise results of~\cite{gim}, obtained using generating functions,
the exact limiting probability for the event that $P_{n}$
will contain any given component is also known.

More recently,
attention has turned to the graph $P_{n,m}$
taken uniformly at random from the set of all labelled planar graphs on $\{1,2,\ldots,n\}$
with exactly $m(n)$ edges,
and the probability that $P_{n,m}$ will contain given components/subgraphs
has been investigated in~\cite{dow2} and~\cite{ger}.

Clearly,
$P_{n,m}$ can be thought of as a random planar graph with a restriction on the average degree.
In this paper,
we shall instead study a random planar graph with restrictions on the maximum and minimum degrees,
again investigating the probability that such a graph will contain given components/subgraphs.
As with $P_{n}$ and $P_{n,m}$,
we shall work solely with labelled graphs.

Given functions $d(n)$ and $D(n)$,
let $\mathcal{P}(n,d,D)$ denote the set of all labelled planar graphs on $\{ 1,2, \ldots, n \}$
with maximum degree at most $D(n)$ and minimum degree at least $d(n)$
(i.e.~with all degrees between $d(n)$ and $D(n)$, inclusive)
and let $P_{n,d,D}$ denote a graph taken uniformly at random from this set
(thus, $P_{n,0,n-1} = P_{n}$, the standard random planar graph).
Note that all graphs with maximum degree at most $2$ are planar,
and so the interest lies with the case when $D(n) \geq 3$ for all $n$.

The structure of this paper shall be based on that of \cite{mcd},
where the standard random planar graph was studied.
Hence, we will start in Section~\ref{bconn}
by establishing a lower bound for the probability that $P_{n,d,D}$ will be connected
(and, hence, an upper bound for the probability that it will contain any given component).
In Section~\ref{growth} we shall use this to show
that there exists a non-zero finite `growth constant' for $|\mathcal{P}(n,d,D)|$,
and in Section~\ref{apps} we will use this second fact to show
that $P_{n,d,D}$ is likely to have many special `appearance'-type copies of certain graphs.
In Section~\ref{cpts}, we shall then use this last result to
deduce a lower bound for the probability that $P_{n,d,D}$ will contain given components.
Finally, in Section~\ref{subs},
we will prove that
$P_{n,d,D}$ has linearly many copies of most, but not all,
planar $H$ satisfying $\Delta(H) \leq \liminf D(n)$
(see Table~\ref{table}).

\begin{table} [ht]
\begin{tabular} {|c|}
\hline
$\phantom{p} \delta (H) = D(n) \textrm{ for all } n
\Rightarrow
\left\{ \begin{array}{l}
\liminf \mathbf{P} > 0 \textrm{ (Theorem~\ref{bounded404})} \\
\& \limsup \mathbf{P} < 1 \textrm{ (Theorem~\ref{bounded311})} \\
\end{array} \right.$ \\
\\
$\left.\begin{array}{l}
\delta (H) < D(n) \textrm{ for all } n \\
\& \not\exists
n:d(n)=D(n)=4
\end{array} \right\}
\Rightarrow
\mathbf{P} \to 1 \textrm{ (Theorem~\ref{bounded1001})}$ \\
\\
$\left. \begin{array}{l}
\delta (H) < D(n) \textrm{ for all } n \\
\& \textrm{ } d(n) = D(n) = 4 \textrm{ for all } n
\end{array} \right\}
\Rightarrow
\left\{ \begin{array}{l}
\mathbf{P} \to 1 \textrm{ if $\exists$ $4$-regular planar }
G \supset H \textrm{ (Theorem\ref{bounded1001})} \\
\mathbf{P}=0 \textrm{ otherwise (trivial observation)} \\
\end{array} \right.$ \\
\hline
\end{tabular}
\caption{A description of
$\mathbf{P} := \mathbf{P}[P_{n,d,D}$ will have a copy of $H$]
for connected planar $H$ with $\Delta (H) \leq \liminf D(n)$.} \label{table}
\end{table}

\section{Connectivity} \label{bconn}

We will start by examining the probability that our random graph is connected.
Not only is this topic interesting in its own right,
but the results given here will also be important ingredients in later sections.

Recall that we must have $D(n) \geq 3$ for planarity to have any impact.
The main result of this section (Theorem~\ref{bounded311}) will be to show that,
given this,
the probability that $P_{n,d,D}$ will be connected
is bounded away from $0$.

The proof will be based on that of Theorem 2.1 of~\cite{mcd},
but will be slightly more complicated,
as the bound on the maximum degree means that $\mathcal{P}(n,d,D)$ is not edge-addable
(i.e.~the class $\mathcal{P}(n,d,D)$
is not closed under the operation of inserting an edge between two components).
Hence, we shall first prove a very helpful result on short cycles:

\begin{Lemma} \label{bounded1}
Let $k < \frac{1}{15}$,
and let $S$ be a planar graph with at most $k|S|$ vertices of degree $\leq 2$.
Then $S$ must contain at least $\left( \frac{1-15k}{28} \right) |S|$ cycles of size $\leq 6$.
In particular, if $S$ has at most $\frac{|S|}{43}$ vertices of degree $\leq 2$
then $S$ must contain at least $\frac{|S|}{43}$ cycles of size $\leq 6$.
\end{Lemma}
\textbf{Proof}
Fix a planar embedding of $S$.
We shall first show that this embedding must have at least
$\left( \frac{1-15k}{14} \right) |S|$ \textit{faces} of size $\leq 6$
(where, as usual, the `size' of a face denotes the number of edges in the associated facial boundary,
with an edge counted twice if it appears twice in the boundary),
and we will later deduce the lemma from this fact.

We shall argue by contradiction.
Let $f_{i}$ denote the number of faces of size $i$ and
suppose that $\sum_{i \leq 6} f_{i} < \left( \frac{1-15k}{14} \right) |S|$.
We have
\begin{eqnarray*}
2e(S) & = & \sum_{i} i f_{i} \\
& \geq & 7 \sum_{i \geq 7} f_{i} \\
& > & 7 \left( \sum_{i} f_{i} - \left( \frac{1-15k}{14} \right) |S| \right),
\textrm{ by our supposition} \\
& = & 7 \left( e(S) - |S| + \kappa(S) + 1- \left( \frac{1-15k}{14} \right) |S| \right),
\textrm{ by Euler's formula} \\
& > & 7 \left( e(S) - \left( \frac{15(1-k)}{14} \right) |S| \right).
\end{eqnarray*}
Thus,
$\left( \frac{15(1-k)}{2} \right) |S| > 5e(S)$.
But $e(S) \geq \frac{3(1-k)|S|}{2}$,
since $S$ contains at least $(1-k)|S|$ vertices of degree $\geq 3$,
and so $5e(S) \geq \left( \frac{15(1-k)|S|}{2} \right)$.
Thus, we obtain our desired contradiction,
and so it must be that we have at least $\left( \frac{1-15k}{14} \right) |S|$ faces of size $\leq 6$.

Let us now consider these faces of size $\leq 6$.
Note that the boundary of a face of size $\leq 6$ must contain a cycle of size $\leq 6$
as a subgraph unless it is acyclic,
in which case it must be the entire graph $S$.
But if $S$ were acyclic,
then at least half of the vertices would have degree $\leq 2$
(since we would have $e(S) \leq |S|-1$),
and this would contradict the conditions of this lemma.
Thus, for each of our faces of size $\leq 6$,
the boundary must contain a cycle of size $\leq 6$ as a subgraph.

Each edge of $S$ can only be in at most two faces of the embedding,
and so each cycle can only be in at most two faces.
Thus, $S$ must contain at least $\left( \frac{1-15k}{28} \right) |S|$ \textit{distinct} cycles of size $\leq 6$.
$\phantom{qwerty}$
\setlength{\unitlength}{0.25cm}
\begin{picture}(1,1)
\put(0,0){\line(1,0){1}}
\put(0,0){\line(0,1){1}}
\put(1,1){\line(-1,0){1}}
\put(1,1){\line(0,-1){1}}
\end{picture} \\
\\

We shall now use Lemma~\ref{bounded1} to obtain our aforementioned main result:

\begin{Theorem} \label{bounded311}
There exists a constant $c >0$ such that
\begin{displaymath}
\mathbf{P} [P_{r,d,D} \textrm{ will be connected}] > c
\end{displaymath}
for all constants $r,d,D \in \mathbf{N} \cup \{ 0 \}$ satisfying $D \geq 3$
and $\mathcal{P}(r,d,D) \neq \emptyset$.
\end{Theorem}
\textbf{Sketch of Proof}
We shall show that there are many ways to construct a graph in $\mathcal{P}(r,d,D)$ with $k-1$ components
from a graph in $\mathcal{P}(r,d,D)$ with $k$ components,
by combining two components.
Our stated lower bound for the proportion of graphs with exactly one component will then follow
by `cascading' this result downwards.

If $D > 6$,
we shall see that we may obtain sufficiently many ways to combine components
simply by inserting edges between them that don't interfere with this upper bound on the maximum degree.

If $D \leq 6$,
we will sometimes also delete an edge from a small cycle
in order to maintain $\Delta \leq D$.
We shall use Lemma~\ref{bounded1} to show that
we have lots of choices for this small cycle,
and then the fact that it is small
(combined with the knowledge that $D<7$)
will help us to bound the amount of double-counting. \\
\\
\textbf{Full Proof}
Choose any $r,d,D \in \mathbf{N} \cup \{ 0 \}$ with $D \geq 3$ and $\mathcal{P}(r,d,D) \neq \emptyset$.
We shall show that there exists a strictly positive constant $c$,
independent of $r,d$ and $D$,
such that
$\mathbf{P} [P_{r,d,D} \textrm{ will be connected}] > c$.

Let $\mathcal{P}^{t}(r,d,D)$ denote the set of graphs in
$\mathcal{P}(r,d,D)$ with exactly $t$ components.
For $k>1$, we shall construct graphs in $\mathcal{P}^{k-1}(r,d,D)$
from graphs in $\mathcal{P}^{k}(r,d,D)$.

Let the graph $G \in \mathcal{P}^{k}(r,d,D)$
and let us denote the $k$ components of $G$ by $S_{1}, S_{2}, \ldots,S_{k}$,
where $|S_{i}|=n_{i}$ for all $i$.
Without loss of generality,
we may assume that $S_{1}, S_{2}, \ldots, S_{k}$ are ordered so that
$S_{i}$ contains at least $\frac{n_{i}}{43}$ vertices of degree $<D$ iff $i \leq l$,
for some fixed $l \in \{ 0,1, \ldots, k \}$.
Note that we must have $l=k$ if $D > 6$,
since (by planarity) $e(S_{i}) < 3n_{i}$
and so we can only have at most $\frac{6n_{i}}{7}$ vertices of degree $\geq 7$.

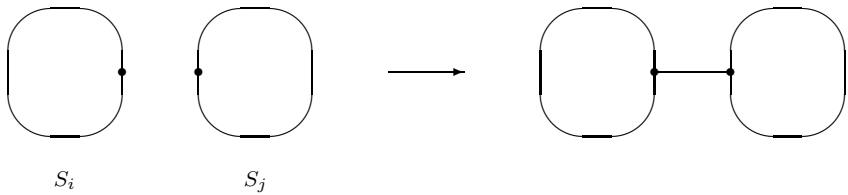
\begin{figure} [ht]
\setlength{\unitlength}{1cm}
\begin{picture}(20,2.35)(-1,-0.75)

\put(-0.5,0.5){\line(0,1){0.5}}
\put(1,0.5){\line(0,1){0.5}}
\put(2,0.5){\line(0,1){0.5}}
\put(3.5,0.5){\line(0,1){0.5}}
\put(6.5,0.5){\line(0,1){0.5}}
\put(8,0.5){\line(0,1){0.5}}
\put(9,0.5){\line(0,1){0.5}}
\put(10.5,0.5){\line(0,1){0.5}}

\put(0.25,1){\oval(1.5,1.2)[t]}
\put(2.75,1){\oval(1.5,1.2)[t]}
\put(7.25,1){\oval(1.5,1.2)[t]}
\put(9.75,1){\oval(1.5,1.2)[t]}

\put(0.25,0.5){\oval(1.5,1.2)[b]}
\put(2.75,0.5){\oval(1.5,1.2)[b]}
\put(7.25,0.5){\oval(1.5,1.2)[b]}
\put(9.75,0.5){\oval(1.5,1.2)[b]}

\put(4.5,0.75){\vector(1,0){1}}

\put(8,0.75){\line(1,0){1}}

\put(1,0.75){\circle*{0.1}}
\put(2,0.75){\circle*{0.1}}
\put(8,0.75){\circle*{0.1}}
\put(9,0.75){\circle*{0.1}}

\put(0.1,-0.75){$S_{i}$}
\put(2.6,-0.75){$S_{j}$}

\end{picture}

\caption{Constructing the graph $G_{i,j}$ in case (a).} \label{bconnfig}
\end{figure}

For $1 \leq i < j \leq k$,
let us construct a new graph $G_{i,j} \in \mathcal{P}^{k-1}(r,d,D)$ as follows: \\
\\
Case (a): if $j \leq l$ (note that this is always the case if $D > 6$) \\
Insert an edge between a vertex in $S_{i}$ of degree $<D$
(we have at least $\frac{n_{i}}{43}$ choices for this)
and a vertex in $S_{j}$ of degree $<D$
(we have at least $\frac{n_{j}}{43}$ choices for this).
See Figure~\ref{bconnfig}.

\begin{figure} [ht]
\setlength{\unitlength}{1cm}
\begin{picture}(20,2.35)(-1,-0.75)

\put(-0.5,0.5){\line(0,1){0.5}}
\put(1,0.5){\line(0,1){0.5}}
\put(2,0.5){\line(0,1){0.5}}
\put(3.5,0.5){\line(0,1){0.5}}
\put(6.5,0.5){\line(0,1){0.5}}
\put(10.5,0.5){\line(0,1){0.5}}

\put(0.25,1){\oval(1.5,1.2)[t]}
\put(2.75,1){\oval(1.5,1.2)[t]}
\put(7.25,1){\oval(1.5,1.2)[t]}
\put(9.75,1){\oval(1.5,1.2)[t]}

\put(0.25,0.5){\oval(1.5,1.2)[b]}
\put(2.75,0.5){\oval(1.5,1.2)[b]}
\put(7.25,0.5){\oval(1.5,1.2)[b]}
\put(9.75,0.5){\oval(1.5,1.2)[b]}

\put(4.5,0.75){\vector(1,0){1}}

\put(8,0.5){\line(1,0){1}}
\put(8,1){\line(1,0){1}}

\put(1,0.5){\circle*{0.1}}
\put(2,0.5){\circle*{0.1}}
\put(8,0.5){\circle*{0.1}}
\put(9,0.5){\circle*{0.1}}
\put(1,1){\circle*{0.1}}
\put(2,1){\circle*{0.1}}
\put(8,1){\circle*{0.1}}
\put(9,1){\circle*{0.1}}

\put(2,0.75){\oval(1,0.5)[r]}
\put(9,0.75){\oval(1,0.5)[r]}

\put(1.75,1.1){$u$}
\put(1.75,0.2){$v$}
\put(1.1,1.1){$w$}
\put(1.1,0.2){$x$}
\put(8.75,1.1){$u$}
\put(8.75,0.2){$v$}
\put(8.1,1.1){$w$}
\put(8.1,0.2){$x$}

\put(0.1,-0.75){$S_{i}$}
\put(2.6,-0.75){$S_{j}$}

\end{picture}

\caption{Constructing the graph $G_{i,j}$ in case (b).} \label{connfigb}
\end{figure}
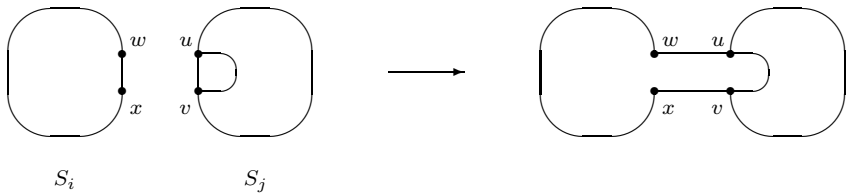

The constructed graph $G_{i,j}$ (see Figure~\ref{bconnfig}) is planar and has exactly $k-1$ components.
It is also clear that we have $d \leq \delta (G_{i,j}) \leq \Delta (G_{i,j}) \leq D$,
since we have not deleted any edges from the original graph
and have only inserted an edge between two vertices with degree $<D$.
Thus, $G_{i,j} \in \mathcal{P}^{k-1}(r,d,D)$. \\
\\
Case (b): if $j>l$ and $n_{i}>1$ (in which case $D \leq 6$) \\
If $j>l$,
then $S_{j}$ contains less than $\frac{n_{j}}{43}$ vertices of degree $<D$.
Thus, by Lemma~\ref{bounded1},
$S_{j}$ must contain at least $\frac{n_{j}}{43}$ cycles of size $\leq 6$.
Delete an edge $uv$ in one of these cycles
(we have at least
$\frac{3}{D+D^{2}+D^{3}+D^{4}}\frac{n_{j}}{43} \geq \frac{3}{6+6^{2}+6^{3}+6^{4}}\frac{n_{j}}{43}$
choices for this edge,
since each cycle must contain at least
$3$ edges and each edge is in at most $(D-1)^{m-2}<D^{m-2}$ cycles of size $m$),
insert an edge between $u$ and a vertex $w \in S_{i}$
(we have $n_{i}$ choices for $w$),
delete an edge between $w$ and $x \in \Gamma(w)$
(we have at least one choice for $x$, since $n_{i}>1$),
and insert an edge between $x$ and $v$
(planarity is preserved,
since we may draw $S_{j}$ so that the face containing $u$ and $v$ is on the outside,
and similarly we may draw $S_{i}$ so that the face containing $w$ and $x$ is on the outside).
See Figure~\ref{connfigb}.

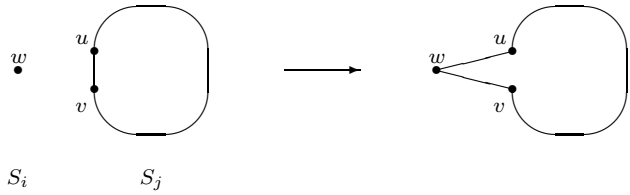
\begin{figure} [ht]
\setlength{\unitlength}{1cm}
\begin{picture}(20,2.35)(-1,-0.75)

\put(2,0.5){\line(0,1){0.5}}
\put(3.5,0.5){\line(0,1){0.5}}
\put(9,0.5){\line(0,1){0.5}}

\put(2.75,1){\oval(1.5,1.2)[t]}
\put(8.25,1){\oval(1.5,1.2)[t]}

\put(2.75,0.5){\oval(1.5,1.2)[b]}
\put(8.25,0.5){\oval(1.5,1.2)[b]}

\put(4.5,0.75){\vector(1,0){1}}

\put(6.5,0.75){\line(4,1){1}}
\put(6.5,0.75){\line(4,-1){1}}

\put(1,0.75){\circle*{0.1}}
\put(2,0.5){\circle*{0.1}}
\put(2,1){\circle*{0.1}}
\put(6.5,0.75){\circle*{0.1}}
\put(7.5,0.5){\circle*{0.1}}
\put(7.5,1){\circle*{0.1}}

\put(1.75,1.1){$u$}
\put(1.75,0.2){$v$}
\put(0.9,0.85){$w$}
\put(7.25,1.1){$u$}
\put(7.25,0.2){$v$}
\put(6.4,0.85){$w$}

\put(0.85,-0.75){$S_{i}$}
\put(2.6,-0.75){$S_{j}$}

\end{picture}

\caption{Constructing the graph $G_{i,j}$ in case (c).} \label{connfigc}
\end{figure}

Since the deleted edge $uv$ was in a cycle,
it was not a cut-edge,
and so the vertex set $V(S_{j})$ is still connected.
The deleted edge $wx$ may have been a cut-edge in $S_{i}$,
but since we have also inserted edges from $w$ to $u \in V(S_{j})$ and from $x$ to $v \in V(S_{j})$
it must be that the vertex set $V(S_{i}) \cup V(S_{j})$ is now connected.
Thus, the constructed planar graph $G_{i,j}$ has exactly $k-1$ components.
By construction, the degrees of the vertices have not changed,
and so we have $d \leq \delta (G_{i,j}) \leq \Delta (G_{i,j}) \leq D$.
Thus, $G_{i,j} \in \mathcal{P}^{k-1}(r,d,D)$. \\
\\
Case (c): if $j>l$ and $n_{i}=1$ (in which case $D \leq 6$) \\
Delete any edge $uv$ in $S_{j}$
(we have at least $n_{j}$ choices for this,
since $S_{j}$ cannot be a forest if $j>l$)
and insert edges $uw$ and $vw$,
where $w$ is the unique vertex in $S_{i}$
(see Figure~\ref{connfigc}).

The constructed graph $G_{i,j}$ is planar and has exactly $k-1$ components.
Note that the degrees have not changed,
except that we now have deg$(w)=2$.
But since $D \geq 3$,
we still have $d \leq \delta(G_{i,j}) \leq \Delta(G_{i,j}) \leq D$.
Thus, we have $G_{i,j} \in \mathcal{P}^{k-1}(r,d,D)$.
\\

Let $z= \frac{3}{43(6+6^{2}+6^{3}+6^{4})}
= \min \left\{ \left( \frac{1}{43} \right)^{2}, \frac{3}{43(6+6^{2}+6^{3}+6^{4})} \right\}$.
Then in all cases
we have at least $zn_{i}n_{j}$ choices when constructing the new graph.
Thus, from our initial graph $G$,
we have at least $\sum_{i<j} zn_{i}n_{j} = z \sum_{i<j} n_{i}n_{j}$
ways to construct a graph in $\mathcal{P}^{k-1}(r,d,D)$.
Note that if $x \leq y$ then $xy > (x- 1)(y+1)$,
so $\sum_{i<j} n_{i}n_{j}$
is at least what it would be if one component in $G$ had order $r-(k-1)$
and the other $k-1$ components were isolated vertices.
Thus,
$z \sum_{i<j} n_{i}n_{j} \geq z \left( \frac{1}{2}(k-1)(k-2) + (k-1)(r-k+1) \right)
= (k-1) \left( r-\frac{k}{2} \right)z$.
Hence, for $k> 1$,
we have at least
$(k-1) \left( r-\frac{k}{2} \right) z|\mathcal{P}^{k}(r,d,D)|
\geq (k-1) \frac{r}{2} z|\mathcal{P}^{k}(r,d,D)|$
ways to construct (not necessarily distinct) graphs in $\mathcal{P}^{k-1}(r,d,D)$. \\

Let us now consider the amount of double-counting:

Given one of our constructed graphs,
there are at most $3$ possibilities for how the graph was obtained (case (a), (b) or (c)).

If case (a) was used (which must be so if $D > 6$),
then we can re-obtain the original graph simply by deleting the inserted edge,
for which there are at most $r-(k-1)<r$ possibilities,
since it must now be a cut-edge.
Thus, if case (a) was used,
we have less than $r$ possibilities for the original graph.

If case (b) was used,
then we can re-obtain the original graph by locating the vertices $u,v,w$ and $x$,
deleting the two inserted edges ($uw$ and $vx$)
and re-inserting the two deleted edges ($uv$ and $wx$).
Note that we have at most $r$ possibilities for which vertex is $u$.
We know that $u$ and $v$ were originally on a cycle of size $\leq 6$,
and so $v$ is still at distance at most $5$ from $u$.
Since the graph has maximum degree at most $D$,
we therefore have at most $D^{2}+D^{3}+D^{4}+D^{5}$ possibilities for $v$.
Once we have located $u$ and $v$,
we then have at most $D$ possibilities for $w$
and at most $D$ possibilities for $x$,
since $w$ and $x$ are now neighbours of $u$ and $v$, respectively.
Thus, if case (b) was used,
we have at most
$D^{2}(D^{2}+D^{3}+D^{4}+D^{5})r \leq 36(6^{2}+6^{3}+6^{4}+6^{5})r$
possibilities for the original graph.

If case (c) was used,
then we can re-obtain the original graph by locating the vertices $u,v$ and $w$,
deleting the two inserted edges ($uw$ and $vw$)
and re-inserting the deleted edge ($uv$).
We have at most $r$ possibilities for which vertex is $w$,
and given $w$ we then know which edges to delete and insert,
as $v$ and $w$ are the only vertices adjacent to $u$.
Thus, if case (c) was used,
we have at most $r$ possibilities for the original graph.

Therefore,
there are less than $r$ possibilities for the original graph if $D > 6$,
since case (a) must have been used,
and less than
$r + 36(6^{2}+6^{3}+6^{4}+6^{5})r +r = 2r(1+18(6^{2}+6^{3}+6^{4}+6^{5}))$
possibilities for the original graph if $D \leq 6$,
since any of case (a), case (b) or case (c) may have been used. \\

Let
$\alpha = \frac{z}{4(1+18(6^{2}+6^{3}+6^{4}+6^{5}))}
= \min \left\{ \frac{z}{2}, \frac{z}{4(1+18(6^{2}+6^{3}+6^{4}+6^{5}))} \right\}$.
Then we have shown that
we can construct at least
$\alpha(k-1) |\mathcal{P}^{k}(r,d,D)|$
\textit{distinct} graphs in $\mathcal{P}^{k-1}(r,d,D)|$,
and so
$
|\mathcal{P}^{k-1}(r,d,D)| \geq \alpha(k-1) |\mathcal{P}^{k}(r,d,D)| \textrm{ for all } k>1.
$

Let us define
$p_{k}$
to be
$\frac{|\mathcal{P}^{k+1}(r,d,D)|}{|\mathcal{P}(r,d,D)|}$
and let
$p=p_{0}= \frac{|\mathcal{P}^{1}(r,d,D)|}{|\mathcal{P}(r,d,D)|} =
\mathbf{P}[P_{r,d,D}$ will be connected].
From the previous paragraph,
we know 
$|\mathcal{P}^{k+1}(r,d,D)| \leq \frac{|\mathcal{P}^{k}(r,d,D)|}{\alpha k}$
for all $k>0$,
and so $p_{k} \leq \frac{p}{\alpha^{k}k!}$ for all $k \geq 0$.
We must have $\sum_{k \geq 0} p_{k} = 1$,
so $\sum_{k \geq 0} \frac{p}{\alpha^{k}k!} \geq 1$
and hence
$p \geq \left( \sum_{k \geq 0} \frac{\left( \frac{1}{\alpha} \right)^{k}}{k!} \right)^{-1}
= e^{-\frac{1}{\alpha}}$.~
\phantom{}
\setlength{\unitlength}{0.25cm}
\begin{picture}(1,1)
\put(0,0){\line(1,0){1}}
\put(0,0){\line(0,1){1}}
\put(1,1){\line(-1,0){1}}
\put(1,1){\line(0,-1){1}}
\end{picture} \\
\\

\section{Growth Constants} \label{growth}

We shall now look at the topic of `growth constants',
which will play a vital role in the proofs of Section~\ref{apps}.

It is known from~\cite{mcd} that there exists a finite constant $\gamma_{l} > 0$
such that $\left( \frac{|\mathcal{P}(n,0,n-1)|}{n!} \right)^{1/n} \to \gamma_{l}$ as $n \to \infty$.
In this section,
we shall use our connectivity bound from Theorem~\ref{bounded311}
to also obtain (in Theorems~\ref{bounded104} and~\ref{bounded105})
growth constants for $\mathcal{P}(n,d,D)$
for the case when $d(n)$ is a constant and $D(n)$ is any monotonically non-decreasing function
(it will turn out that the result for this restricted case is all that will be required for later sections).

We shall follow the proof of Theorem 3.3 of~\cite{mcd},
which will require us to first state the following useful lemma:

\begin{Proposition}[see, for example, Lemma 11.6 of \cite{vanL}] \label{vanL 11.6}
Let $f:\mathbf{N} \to \mathbf{R}^{+}$ be a function such that
$f(n)>0$ for all large $n$ and
$f(i+j) \geq f(i) \cdot f(j)$ for all $i,j \in \mathbf{N}$.
Then $(f(n))^{1/n} \to \sup_{n} \left( (f(n))^{1/n} \right)$ as $n \to \infty$.
\end{Proposition}

We may now use Proposition~\ref{vanL 11.6} to obtain our growth constant result:

\begin{Theorem} \label{bounded104}
Let $d \in \{ 0,1, \ldots, 5 \}$ be a constant and
let $D(n)$ be a monotonically non-decreasing integer-valued function that for all large $n$ satisfies
$D(n) \geq \max \{ d, 3 \}$ and $(d,D(n)) \notin \{ (3,3), (5,5) \}$.
Then there exists a \emph{finite} constant $\gamma_{d,D} >0$ such that
\begin{displaymath}
\left( \frac{|\mathcal{P}(n,d,D)|}{n!} \right)^{\frac{1}{n}} \to
\gamma_{d,D}
\textrm{ as } n \to \infty.
\end{displaymath}
\end{Theorem}
\textbf{Proof}
We shall follow the method of proof of Theorem 3.3 of~\cite{mcd}.
Let $c$ be the constant given by Theorem~\ref{bounded311} and
let $g(n,d,D) = \frac{c^{2}|\mathcal{P}(n,d,D)|}{2 \cdot n!}$ for all $n \in \mathbf{N}$.
We shall show that $g(n,d,D)$ satisfies the conditions of Proposition~\ref{vanL 11.6},
which we will then use to deduce our result.

To show $g(n,d,D)>0$ for all large $n$,
it clearly suffices to prove that $\mathcal{P}(n,d,D)$ is non-empty for all large $n$.
This is true, but fairly tedious to demonstrate
and so we will omit the details
(see~\cite{dow} for a full discussion).

Let us now show that $g$ satisfies the supermultiplicative condition:

Let $i,j \in \mathbf{N}$ and
let us denote by
$\mathcal{P}_{c}(i,d,D)$ and $\mathcal{P}_{c}(j,d,D)$
the set of connected graphs in $\mathcal{P}(i,d,D)$ and $\mathcal{P}(j,d,D)$,
respectively.
Then, by Theorem~\ref{bounded311},
we know that there exists a constant $c>0$ such that we have
$|\mathcal{P}_{c}(i,d,D)| \geq c |\mathcal{P}(i,d,D)|$
and
$|\mathcal{P}_{c}(j,d,D)| \geq c |\mathcal{P}(j,d,D)|$.
We may form a graph in $\mathcal{P}(i+j,d,D)$ by choosing $i$ of the $i+j$ vertices
$\left( \left(^{i+j}_{\phantom{q}j} \right) \textrm{ choices} \right)$,
placing a connected planar graph $G_{1}$ with $|G_{1}|=i$ and $d \leq \delta(G_{1}) \leq \Delta(G_{1}) \leq D(i)$
on the chosen vertices
($|\mathcal{P}_{c}(i,d,D)| \geq c |\mathcal{P}(i,d,D)|$ choices),
and then placing a connected planar graph $G_{2}$ with $|G_{2}|=j$ and
$d \leq \delta(G_{2}) \leq \Delta(G_{2}) \leq D(j)$
on the remaining $j$ vertices
($|\mathcal{P}_{c}(j,d,D)| \geq c|\mathcal{P}(j,d,D)|$ choices).
If $i=j$,
then we need to divide by two to avoid double-counting.
Note that the constructed graph will have maximum degree at most $\max \{ D(i), D(j) \}$
and so will indeed be in $\mathcal{P}(i+j,d,D)$,
since $D$ is a monotonically non-decreasing function.
Thus,
\begin{eqnarray*}
|\mathcal{P}(i+j,d,D)| \geq
\frac{c^{2}}{2} \left(^{i+j}_{\phantom{q}j} \right)
|\mathcal{P}(i,d,D)| \cdot |\mathcal{P}(j,d,D)| \textrm{ for all } i,j
\end{eqnarray*}
and, therefore,
\begin{eqnarray*}
g(i+j,d,D) & = & \frac{c^{2} |\mathcal{P}(i+j,d,D)|}{2(i+j)!} \\
& \geq & \frac{ c^{4} \left(^{i+j}_{\phantom{q}j} \right)
|\mathcal{P}(i,d,D)| |\mathcal{P}(j,d,D)|}
{4(i+j)!} \\
& = & \frac{ c^{2} |\mathcal{P}(i,d,D)|}{2 \cdot i!}
\frac{c^{2} |\mathcal{P}(j,d,D)|}{2 \cdot j!} \\
& = & g(i,d,D) \cdot g(j,d,D).
\end{eqnarray*}

Let $\gamma_{d,D} = \sup_{n} \left( (g(n,d,D))^{1/n} \right)$.
By Proposition~\ref{vanL 11.6},
it now only remains to show that $\gamma_{d,D} < \infty$.
But clearly $\mathcal{P}(n,d,D) \subset \mathcal{P}(n,0,n-1)$,
the set of all labelled planar graphs on $\{ 1,2, \ldots, n \}$,
and it is known that
$\left( \frac{|\mathcal{P}(n,0,n-1)|}{n!} \right)^{\frac{1}{n}}$
converges to a finite constant as $n \to \infty$
(see~\cite{mcd}),
so we are done.
\phantom{qwerty}
\setlength{\unitlength}{0.25cm}
\begin{picture}(1,1)
\put(0,0){\line(1,0){1}}
\put(0,0){\line(0,1){1}}
\put(1,1){\line(-1,0){1}}
\put(1,1){\line(0,-1){1}}
\end{picture} \\
\\

By the same proof,
we may also obtain an analogous result to Theorem~\ref{bounded104}
for the case when $D(n) = d \in \{ 3,5 \}$ for all $n$:

\begin{Theorem} \label{bounded105}
Let $D \in \{ 3,5 \}$ be a constant.
Then there exists a finite constant $\gamma_{D,D}>0$ such that
\begin{displaymath}
\left( \frac{|\mathcal{P}(2n,D,D)|}{(2n)!} \right)^{\frac{1}{2n}} \to
\gamma_{D,D} \textrm{ as } n \to \infty.
\end{displaymath} \\
\end{Theorem}

\section{Appearances} \label{apps}

We shall now look at special subgraphs in $P_{n,d,D}$ called `appearances',
with the aim of turning some of these into components in Section~\ref{cpts}.

We will produce separate results for the cases when we have
$d(n) < D(n)$ for all $n$ (Theorem~\ref{bounded11})
and when $d(n) = D(n)$ for all $n$ (Theorem~\ref{bounded110}),
because for the latter case it will be more awkward
to later convert these subgraphs into components
without violating our bound on the minimum degree.

We will deal with the $d(n) < D(n)$ case first.
The main work will be done in Lemma~\ref{bounded106},
where we shall follow the lines of a proof of \cite{mcd}
(using the growth constants of Section~\ref{growth})
to obtain an appearance result for the case when
$d(n)$ is a constant and $D(n)$ is a monotonically non-decreasing function.
We shall then extend this into a more general form in Theorem~\ref{bounded11},
before finally noting that the $d(n)=D(n)$ case follows from a similar proof. \\

We start with a definition of appearances:

\begin{Definition} \label{defapps}
Let $H$ be a graph on the vertex set $\{1,2,\ldots,h\},$ and let $G$ be a graph on the vertex set $\{1,2,\ldots,n\}$,
where $n>h$.
Let $W \subset V(G)$ with $|W|=h$, and let the `root' $r_{W}$ denote the least element in $W$.
We say that $H$ \emph{appears} at $W$ in $G$ if
(a) the increasing bijection from ${1,2,\ldots,h}$ to $W$ gives
an isomorphism between $H$ and the induced subgraph $G[W]$ of $G$;
and (b) there is exactly one edge in $G$ between $W$ and the rest of $G$,
and this edge $e_{W} = r_{W}v_{W}$ is incident with the root $r_{W}$.

We call $e_{W}$ the \emph{associated cut-edge} of the appearance,
and we say that $TE_{W} := E(G[W]) \cup \{e_{W}\}$ is the \emph{total edge set} of the appearance
(see Figure~\ref{cutapp}).

When working in $\mathcal{P}(n,d,D)$,
we say that the appearance is \emph{cut-able} if we have
$\min \{ \deg_{G}r_{W}, \deg_{G}v_{W} \} > d$,
and we let $f_{H}(G)$ denote the number of cut-able appearances of $H$ in $G$
(that is, the number of sets $W \subset V(G)$ such that there is a cut-able appearance of $H$ at $W$).
\begin{figure} [ht]
\setlength{\unitlength}{1cm}
\begin{picture}(20,4.25)(-3.375,-0.75)

\put(0.75,0){\line(1,0){0.5}}
\put(0.75,1){\line(1,0){0.5}}
\put(0.75,2){\line(1,0){0.5}}
\put(0.75,3.5){\line(1,0){0.5}}
\put(4.25,0){\line(1,0){0.5}}
\put(4.25,1){\line(1,0){0.5}}

\put(1,1){\line(0,1){1}}
\put(4.5,1){\line(0,1){1}}

\put(0.75,0.5){\oval(0.5,1)[l]}
\put(0.75,2.75){\oval(2,1.5)[l]}
\put(4.25,0.5){\oval(0.5,1)[l]}
\put(1.25,0.5){\oval(0.5,1)[r]}
\put(1.25,2.75){\oval(2,1.5)[r]}
\put(4.75,0.5){\oval(0.5,1)[r]}

\put(1,1){\circle*{0.1}}
\put(1,2){\circle*{0.1}}
\put(4.5,1){\circle*{0.1}}
\put(4.5,2){\circle*{0.1}}

\put(0.6,2.7){$G \setminus W$}
\put(0.65,0.4){\small{$G[W]$}}
\put(4.15,0.4){\small{$G[W]$}}

\put(0.9,-0.75){\large{$G$}}
\put(4.1,-0.75){\large{$TE_{W}$}}

\put(0.5,1.4){$e_{W}$}
\put(1.1,1.8){$v_{W}$}
\put(1.1,1.1){$r_{W}$}
\put(4,1.4){$e_{W}$}
\put(4.6,1.8){$v_{W}$}
\put(4.6,1.1){$r_{W}$}

\end{picture}

\caption{$\textrm{An appearance at $W$ in $G$ and its total edge set}$.} \label{cutapp}
\end{figure}
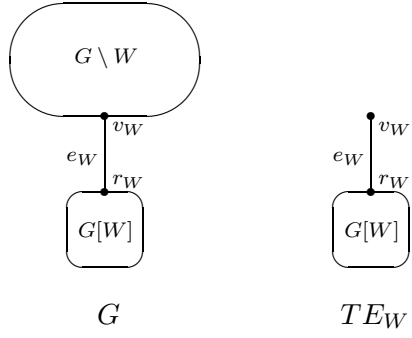
\end{Definition}

We are now ready to give our first result on appearances in $P_{n,d,D}$.
We shall start by assuming that $d(n)$ is constant and $D(n)$ is monotonically non-decreasing,
as in Section~\ref{growth},
but we will later (in Theorem~\ref{bounded11}) get rid of these conditions.
The statement of the result may seem complicated,
but basically it just asserts that for any `sensible' choice of $H$,
there will probably be lots of cut-able appearances of $H$ in $P_{n,d,D}$.
Clearly, `sensible' entails that we must have $\delta(H) \geq d$,
$\Delta(H) \leq D(n)$ and $\deg_{H}(1)+1 \leq D(n)$,
and as always we will require that $D(n) \geq 3$
(note that it follows from these conditions that we must also have $d< D(n)$).
The proof is based on that of Theorem 4.1 of~\cite{mcd}.

\begin{Lemma} \label{bounded106}
Let $H$ be a (fixed) connected planar graph on $\{ 1,2, \ldots, h \}$.
Then there exists a constant $\alpha(h)>0$ such that,
given any constant $d \leq \delta(H)$
and any monotonically non-decreasing integer-valued function $D(n)$ satisfying
$\liminf_{n \to \infty}D(n) \geq \max \{ \Delta(H), \deg_{H}(1)+1, 3 \}$,
we have
\begin{displaymath}
\mathbf{P} [ f_{H}(P_{n,d,D}) \leq \alpha n ] < e^{-\alpha n} \textrm{ for all sufficiently large }n.
\end{displaymath}
\end{Lemma}
\textbf{Sketch of Proof}
We choose a specific $\alpha$
and suppose that the result is false for $n=k$,
where $k$ is suitably large.
Using Theorem~\ref{bounded104},
it then follows that there are many graphs $G \in \mathcal{P}(k,d,D)$ with $f_{H}(G) \leq \alpha k$.

From each such $G$,
we construct graphs in $\mathcal{P}((1+\delta)k,d,D)$, for a fixed $\delta >0$.
If $G$ has lots of vertices with degree $<D(k)$,
then we do this simply by attaching appearances of $H$ to some of these vertices.
If $G$ has few vertices with degree $<D(k)$,
then we attach appearances of $H$ to small cycles in $G$
and also delete appropriate edges.
By Lemma~\ref{bounded1},
we have lots of choices for these small cycles
and, since $G$ has few vertices with degree $<D(k)$,
we may assume that we don't interfere with any vertices of minimum degree.

The fact that the original graphs satisfied $f_{H} \leq \alpha k$,
together with the knowledge that any deleted edges were in small cycles,
is then used to show that there is not much double-counting,
and so we find that we have constructed so many graphs in $\mathcal{P}((1+\delta)k,d,D)$
that we contradict Theorem~\ref{bounded104}. \\
\\
\textbf{Full Proof}
Let $p \in \left( 0, \frac{1}{7(6^{2}+6^{3}+6^{4}+6^{5})} \right)$,
let $\beta = \frac{344e^{2}(h+7)(6^{2}+6^{3}+6^{4}+6^{5})h! \left( \gamma_{d,D} \right)^{h}}{p}$,
and let $\alpha \in
\left( 0,  \frac{p}{344e^{2}(h+7)(6^{2}+6^{3}+6^{4}+6^{5})h! \left( \gamma_{l} \right)^{h}} \right)$,
where we recall that $\gamma_{l} \approx 27.2268$
denotes the growth constant for $\mathcal{P}(n,0,n-1)$.
Clearly $\gamma_{d,D} \leq \gamma_{l}$,
so $\alpha \beta <1$
and hence there exists an
$\epsilon \in \left( 0,\frac{1}{3} \right)$ such that $(\alpha \beta)^{\alpha} = 1-3 \epsilon$.

By Theorem~\ref{bounded104},
there exists an $N$ such that
\begin{equation} \label{eq:gr}
(1-\epsilon)^{n}n! \left( \gamma_{d,D} \right)^{n}
\leq |\mathcal{P}(n,d,D)|
\leq (1+\epsilon)^{n}n! \left( \gamma_{d,D} \right)^{n}
\textrm{ for all } n \geq N.
\end{equation}
Suppose (aiming for a contradiction)
that we can find a value
$k>N$ such that
$\mathbf{P} [ f_{H}(P_{k,d,D}) \leq \alpha k ] \geq e^{-\alpha k}$,
and let $\mathcal{G}$ denote the set of graphs in $\mathcal{P}(k,d,D)$ such that
$G \in \mathcal{G}$ iff $f_{H}(G) \leq \alpha k$.
Then we must have
$|\mathcal{G}| \geq e^{-\alpha k} |\mathcal{P}(k,d,D)| \geq
e^{-\alpha k} (1-\epsilon)^{k}k! \left( \gamma_{d,D} \right)^{k}$.

Let $\delta = \frac{\lceil \alpha k \rceil h}{k}$.
We may assume that $k$ is sufficiently large that
$\lceil \alpha k \rceil \leq 2 \alpha k$.
Thus, $\delta \leq 2 \alpha h < 1$ \label{delta}
(by our definition of $\alpha$).
This fact will be useful later. \\

We shall construct graphs in $\mathcal{P}((1+\delta)k,d,D)$:

Choose $\delta k$ special vertices
(we have $\left(^{(1+\delta)k} _{\phantom{qq} \delta k} \right)$ choices for these)
and partition them into $\lceil \alpha k \rceil$ unordered blocks of size $h$
(we have $\left(^{\phantom{w} \delta k} _{h, \ldots, h} \right) \frac{1}{\lceil \alpha k \rceil !}$ choices for this).
On each of the blocks,
put a copy of $H$ such that the increasing bijection from $\{ 1,2, \ldots, h \}$ to the block
is an isomorphism between $H$ and this copy.
Note that we may assume that $k$ is large enough that $D(k) \geq \liminf_{n \to \infty}D(n)$,
and so the root, $r_{B}$, of a block
(i.e.~the lowest numbered vertex in it)
satisfies deg$(r_{B})<D(k)$,
by the conditions of the theorem.
On the remaining $k$ vertices,
we place a planar graph $G$ with
$d \leq \delta(G) \leq \Delta(G) \leq D(k)$
and $f_{H}(G) \leq \alpha k$
(we have at least $|\mathcal{G}|$ choices for this).

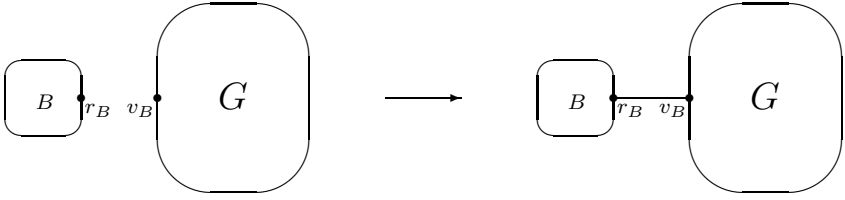
\begin{figure} [ht]
\setlength{\unitlength}{1cm}
\begin{picture}(20,2.5)(-0.5,-0.5)

\put(0,0.5){\line(0,1){0.5}}
\put(1,0.5){\line(0,1){0.5}}
\put(2,0.5){\line(0,1){0.5}}
\put(4,0.5){\line(0,1){0.5}}
\put(7,0.5){\line(0,1){0.5}}
\put(8,0.5){\line(0,1){0.5}}
\put(9,0.5){\line(0,1){0.5}}
\put(11,0.5){\line(0,1){0.5}}

\put(8,0.75){\line(1,0){1}}

\put(5,0.75){\vector(1,0){1}}

\put(0.5,1){\oval(1,0.5)[t]}
\put(0.5,0.5){\oval(1,0.5)[b]}
\put(3,1){\oval(2,2)[t]}
\put(3,0.5){\oval(2,2)[b]}
\put(7.5,1){\oval(1,0.5)[t]}
\put(7.5,0.5){\oval(1,0.5)[b]}
\put(10,1){\oval(2,2)[t]}
\put(10,0.5){\oval(2,2)[b]}

\put(1,0.75){\circle*{0.1}}
\put(2,0.75){\circle*{0.1}}
\put(8,0.75){\circle*{0.1}}
\put(9,0.75){\circle*{0.1}}

\put(0.4,0.6){$B$}
\put(2.8,0.6){\Large{$G$}}
\put(1.05,0.55){\footnotesize{$r_{B}$}}
\put(1.6,0.55){\footnotesize{$v_{B}$}}

\put(7.4,0.6){$B$}
\put(9.8,0.6){\Large{$G$}}
\put(8.05,0.55){\footnotesize{$r_{B}$}}
\put(8.6,0.55){\footnotesize{$v_{B}$}}

\end{picture}

\caption{Creating a cut-able appearance of $H$ at $B$ in case (a).} \label{appa}
\end{figure}

We shall continue our construction in one of two ways,
depending on the number of vertices of degree $D(k)$ in $G$: \\
\\
Case (a): If $G$ has at least $\frac{pk}{43}$ vertices of degree $<D(k)$
(note that this is certainly the case if $D(k) \geq 7$). \\
For each block $B$,
we choose a \textit{different} non-special vertex $v_{B} \in V(G)$ with deg$(v_{B})<D(k)$
(we have at least
$\left( ^{pk/43} _{\phantom{i} \lceil \alpha k \rceil} \right) \lceil \alpha k \rceil !$
choices for this,
since certainly $\alpha < \frac{p}{86}$
and we may assume that $k$ is large enough that
$\lceil \alpha k \rceil \leq 2 \alpha k$),
and we insert the edge $r_{B}v_{B}$ from the root of the block to this vertex,
creating a cut-able appearance of $H$ at $B$ (see Figure~\ref{appa}).
Note that we have not deleted any edges,
so we shall still have minimum degree at least $d$,
and we have only inserted edges between vertices of degree $<D(k)$,
so we still have maximum degree at most $D(k)$,
which is at most $D((1+ \delta)k)$ by monotonicity of $D$.
Thus, our new graph is indeed in $\mathcal{P}((1+\delta)k,d,D)$.

Hence, for each graph $G$ with at least $\frac{pk}{43}$ vertices of degree $<D(k)$,
we find that we can construct at least
$
\left(^{(1+\delta)k} _{\phantom{qq} \delta k} \right)
\left(^{\phantom{p}\delta k}_{h \ldots h} \right)
\cdot \frac{1}{\lceil \alpha k \rceil !} \cdot \left( ^{pk/43}_{\phantom{i} \lceil \alpha k \rceil} \right)
\lceil \alpha k \rceil !
$
different graphs in $\mathcal{P}((1+\delta)k,d,D)$. \\
\\
Case (b): If $G$ has less than $\frac{pk}{43}$ vertices of degree $<D(k)$
(in which case $D(k)<7$). \\
Before describing the case (b) continuation of our construction,
it shall first be useful to investigate the number of short cycles in $G$: \\
If $G$ has less than $\frac{pk}{43} < \frac{k}{43}$ vertices of degree $<D(k)$,
then (by Lemma~\ref{bounded1}) $G$ contains at least $\frac{k}{43}$ cycles of size at most $6$.
A vertex can only be in at most
$(D(k))^{2}+(D(k))^{3}+(D(k))^{4}+(D(k))^{5} \leq 6^{2}+6^{3}+6^{4}+6^{5}$ cycles of size at most $6$,
so $G$ must have at most $\frac{pk(6^{2}+6^{3}+6^{4}+6^{5})}{43}$ cycles of size at most $6$
that contain a vertex of degree $<D(k)$.
In particular, $G$ must have at least $\frac{(1-(6^{2}+6^{3}+6^{4}+6^{5})p)k}{43}$ cycles of size at most $6$
that don't contain a vertex of degree $d$, since $d \leq \delta(H) < \deg_{H}(1)+1 \leq D(k)$.
Since a vertex can only be in at most
$6^{2}+6^{3}+6^{4}+6^{5}$ cycles of size at most $6$,
each cycle of size at most $6$ can only have a vertex in common with
at most $6(6^{2}+6^{3}+6^{4}+6^{5})$ other cycles of size at most $6$.
Thus, $G$ must have a set of at least
$\frac{\left( \frac{1-(6^{2}+6^{3}+6^{4}+6^{5})p}{6(6^{2}+6^{3}+6^{4}+6^{5})} \right)k}{43} > \frac{pk}{43}$
\textit{vertex-disjoint} cycles of size at most $6$ that don't contain a vertex of degree $d$
$\left( \textrm{using the fact that }
p < \frac{1}{7(6^{2}+6^{3}+6^{4}+6^{5})} \right)$.
We shall call these cycles `special'.

\begin{figure} [ht]
\setlength{\unitlength}{1cm}
\begin{picture}(20,2.5)(-0.5,-0.5)

\put(0,0.5){\line(0,1){0.5}}
\put(1,0.5){\line(0,1){0.5}}
\put(2,0.5){\line(0,1){0.5}}
\put(4,0.5){\line(0,1){0.5}}
\put(7,0.5){\line(0,1){0.5}}
\put(8,0.5){\line(0,1){0.5}}
\put(11,0.5){\line(0,1){0.5}}

\put(8,0.75){\line(4,-1){1}}

\put(5,0.75){\vector(1,0){1}}

\put(0.5,1){\oval(1,0.5)[t]}
\put(0.5,0.5){\oval(1,0.5)[b]}
\put(3,1){\oval(2,2)[t]}
\put(3,0.5){\oval(2,2)[b]}
\put(7.5,1){\oval(1,0.5)[t]}
\put(7.5,0.5){\oval(1,0.5)[b]}
\put(10,1){\oval(2,2)[t]}
\put(10,0.5){\oval(2,2)[b]}

\put(2,0.75){\oval(1,0.5)[r]}
\put(9,0.75){\oval(1,0.5)[r]}

\put(1,0.75){\circle*{0.1}}
\put(2,0.5){\circle*{0.1}}
\put(2,1){\circle*{0.1}}
\put(8,0.75){\circle*{0.1}}
\put(9,0.5){\circle*{0.1}}
\put(9,1){\circle*{0.1}}

\put(0.4,0.6){$B$}
\put(3,0.6){\Large{$G$}}
\put(1.05,0.5){\footnotesize{$r_{B}$}}
\put(1.6,0.3){\footnotesize{$v_{B}$}}
\put(1.6,1.15){\footnotesize{$u_{B}$}}

\put(7.4,0.6){$B$}
\put(10,0.6){\Large{$G$}}
\put(8.05,0.5){\footnotesize{$r_{B}$}}
\put(8.6,0.35){\footnotesize{$v_{B}$}}
\put(8.6,1.15){\footnotesize{$u_{B}$}}
\end{picture}

\caption{Creating a cut-able appearance of $H$ at $B$ in case (b).} \label{appb}
\end{figure}
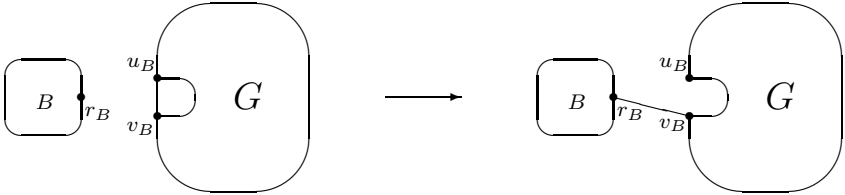

Recall that we have $\lceil \alpha k \rceil$ blocks isomorphic to $H$.
For each block $B$,
choose a \textit{different} one of our `special' cycles
(we have at least $\left( ^{pk/43}_{\phantom{i}\lceil \alpha k \rceil} \right) \lceil \alpha k \rceil !$
choices for this),
delete an edge $u_{B}v_{B}$ in the cycle
and insert an edge $r_{B}v_{B}$ from the root of the block to a vertex $v_{B}$ that was incident to the deleted edge,
creating an appearance of $H$ at $B$ (see Figure~\ref{appb}).
Note that the deleted edge was between two vertices of degree $>d$,
so we still have minimum degree at least $d$
(and $v_{B}$ will still have degree $>d$,
so the appearance will be cut-able).
Recall that the root of each block has degree $<D(k)$,
so we still have maximum degree at most $D(k)$,
which is at most $D((1+ \delta)k)$ by monotonicity of $D$.
Thus, our constructed graph is indeed in $\mathcal{P}((1+\delta)k,d,D)$.

Thus, for each graph $G$ with less than $\frac{pk}{43}$ vertices of degree $<D(k)$,
we find that we can construct at least
$
\left(^{(1+\delta)k} _{\phantom{qq} \delta k} \right)
\left(^{\phantom{p}\delta k}_{h \ldots h} \right)
\cdot \frac{1}{\lceil \alpha k \rceil !} \cdot \left( ^{pk/43}_{\phantom{i} \lceil \alpha k \rceil} \right)
\lceil \alpha k \rceil !
$
different graphs in $\mathcal{P}((1+\delta)k,d,D)$. \\

We have shown that,
regardless of whether case (a) or case (b) is used,
for each $G$ we can construct at least
\begin{eqnarray*}
& & \left(^{(1+\delta)k} _{\phantom{qq} \delta k} \right)
\left(^{\phantom{p}\delta k}_{h \ldots h} \right)
\cdot \frac{1}{\lceil \alpha k \rceil !} \cdot \left( ^{pk/43}_{\phantom{i} \lceil \alpha k \rceil} \right)
\lceil \alpha k \rceil ! \\
& = & \frac{((1+\delta)k)!}{k!} \frac{1}{(h!)^{\lceil \alpha k \rceil}}
\left(^{pk/43}_{\phantom{i} \lceil \alpha k \rceil} \right) \\
& \geq & \frac{((1+\delta)k)!}{k!} \frac{1}{(h!)^{\lceil \alpha k \rceil}}
\frac{(pk/43 - \lceil \alpha k \rceil +1)^{\lceil \alpha k \rceil}}{\lceil \alpha k \rceil !} \\
& \geq & \frac{((1+\delta)k)!}{k!} \frac{1}{(h!)^{\lceil \alpha k \rceil}}
\left( \frac{pk}{86} \right)^{\lceil \alpha k \rceil} \frac{1}{\lceil \alpha k \rceil !} \\
& & \left( \textrm{ since certainly } \alpha < \frac{p}{86}
\textrm{ and so } \frac{pk}{43} - \alpha k \geq \frac{pk}{86} \right) \\
& \geq & \frac{((1+\delta)k)!}{k!} \left( \frac{pk}{86h! \lceil \alpha k \rceil} \right)^{\lceil \alpha k \rceil} \\
& \geq & \frac{((1+\delta)k)!}{k!} \left( \frac{p}{172h! \alpha} \right)^{\lceil \alpha k \rceil} \\
& & \left( \textrm{since we may assume $k$ is large enough that }
\lceil \alpha k \rceil \leq 2 \alpha k \right)
\end{eqnarray*}
different graphs in $\mathcal{P}((1+\delta)k,d,D)$.
Thus, recalling that we have at least
$e^{-\alpha k} (1-\epsilon)^{k}k! \left( \gamma_{d,D} \right)^{k}$
choices for $G$,
we find that we can in total construct at least
$
e^{-\alpha k} (1-\epsilon)^{k} ((1+\delta)k)! \left( \gamma_{d,D} \right)^{k}
\left( \frac{p}{172h! \alpha} \right)^{\lceil \alpha k \rceil}
$
(not necessarily distinct) graphs in $\mathcal{P}((1+\delta)k,d,D)$. \\

We are now at the half way point of our proof,
and it remains to investigate the amount of double-counting,
i.e.~how many times each of our constructed graphs will have been built.
Given one of our constructed graphs, $G^{\prime}$,
there are at most two possibilities for how the graph was obtained (case (a) or case(b)),
and we shall now examine these two cases separately: \\

If case (a) was used, then we can re-obtain the original graph, $G$,
simply by deleting the $\lceil \alpha k \rceil$ cut-able appearances that were deliberately added.
Thus, in order to bound the amount of double-counting
under case (a),
we only need to investigate $f_{H}(G^{\prime})$:

Suppose $W$ is a cut-able appearance of $H$ in $G^{\prime}$.
We shall consider how many possibilities there are for $W$:

(i) If we don't have $TE_{W} \subset E(G)$,
then the total edge set of $W$ must intersect the total edge set of one of our deliberately created appearances.
Note that the total edge set of an appearance of $H$
can only intersect at most $|H|$ other total edge sets of appearances of $H$
(since there are at most $|H|$ cut-edges in the total edge set
and each of these can have at most one `orientation' that provides an appearance of $H$),
so we have at most $(h+1) \lceil \alpha k \rceil$ possibilities for $W$
(including the possibility that $W$ is one of our deliberately created appearances).

If $TE_{W} \subset E(G)$,
then $W$ must have been an appearance of $H$ in $G$:

(ii) If $W$ was a \textit{cut-able} appearance of $H$ in $G$,
then there are at most $\lceil \alpha k \rceil$ possibilities for $W$,
by definition of $\mathcal{G}$.

(iii) If $W$ was an appearance of $H$ in $G$
that was not cut-able,
then the unique vertex $v \in V(G^{\prime}) \setminus W$
incident to the root of $W$ must have had deg$(v) = d$ originally
and must have been chosen as $v_{B}$ by some block $B$.
Hence, we have at most $\lceil \alpha k \rceil$ possibilities for $v$
and thus at most $d \lceil \alpha k \rceil$ possibilities for $W$.

Thus, if case (a) was used,
then $f_{H}(G^{\prime}) \leq (h+d+2) \lceil \alpha k \rceil$,
giving us at most
$\left( ^{(h+d+2) \lceil \alpha k \rceil}_{\phantom{www} \lceil \alpha k \rceil} \right)
\leq ((h+d+2)e)^{\lceil \alpha k \rceil}
\leq ((h+7)e)^{\lceil \alpha k \rceil}$
possibilities for $G$. \\

If case (b) was used,
we can re-obtain the original graph, $G$,
by deleting the $\lceil \alpha k \rceil$ appearances that were deliberately added
and re-inserting the $\lceil \alpha k \rceil$ deleted edges.
Note that once we have identified the appearances that were deliberately added,
we have at most
$\left( (D(k))^{2}+\!(D(k))^{3}+\!(D(k))^{4}+\!(D(k))^{5} \right)^{\lceil \alpha k \rceil}$
$\leq (6^{2}+6^{3}+6^{4}+6^{5})^{\lceil \alpha k \rceil}$
possibilities for the edges that were deleted,
since for each appearance we will automatically know one endpoint, $v$, of the corresponding deleted edge
and we know that the other endpoint, $u$,
will now be at most distance $5$ from $v$,
since $uv$ was originally part of a cycle of size $\leq 6$.
Hence, as with case (a),
it now remains to examine how many possibilities there are for the $\lceil \alpha k \rceil$ appearances
that were deliberately added.

Suppose $W$ is a cut-able appearance of $H$ in $G^{\prime}$.

(i) If we don't have $TE_{W} \subset E(G)$,
then we have at most $(h+1) \lceil \alpha k \rceil$ possibilities for $W$,
as with case (a).

(ii) If $TE_{W} \subset E(G)$
and $W$ was an appearance of $H$ in $G$,
then note that this appearance must have already been cut-able,
since it is clear that we have $\deg_{G^{\prime}} \leq \deg_{G}$
for all vertices that were in $V(G)$.
Hence, there are at most $\lceil \alpha k \rceil$ possibilities for $W$,
by definition of $\mathcal{G}$.

\begin{figure} [ht]
\setlength{\unitlength}{1cm}
\begin{picture}(20,3.5)(-5,0)

\put(0.75,0){\line(1,0){0.5}}
\put(0.75,1){\line(1,0){0.5}}
\put(0.75,2){\line(1,0){0.5}}
\put(0.75,3.5){\line(1,0){0.5}}

\put(1,1){\line(0,1){1}}

\put(0.75,0.5){\oval(0.5,1)[l]}
\put(0.75,2.75){\oval(2,1.5)[l]}
\put(1.25,0.5){\oval(0.5,1)[r]}
\put(1.25,2.75){\oval(2,1.5)[r]}

\put(1.5,0.5){\circle*{0.1}}
\put(1.2,2){\circle*{0.1}}

\put(0.8,0.4){$W$}

\put(0.5,1.4){$e_{W}$}
\put(1.3,1.8){$v_{B}$}
\put(1.6,0.4){$u_{B}$}

\end{picture}

\caption{The appearance $W$ in case (iii).} \label{appiii}
\end{figure}
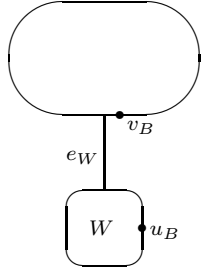

(iii) If $TE_{W} \subset E(G)$
and $W$ was \textit{not} an appearance of $H$ in $G$,
then there must have originally been either another edge between $W$ and $V(G) \setminus W$ other than $e_{W}$,
or another edge between vertices in $W$.
This deleted edge must be of the form $u_{B}v_{B}$ for some block $B$,
and so $W$ must contain either $u_{B}$ or $v_{B}$ (or both).
However, if $v_{B} \in W$ then $r_{B}v_{B}$ would belong to the total edge set of $W$,
which would contradict our assumption that $TE_{W} \subset E(G)$.
Thus, $u_{B} \in W$ and $v_{B} \notin W$ (see Figure~\ref{appiii}).
Recall that the deleted edge $u_{B}v_{B}$ was originally part of a cycle of size $\leq 6$
and that no other edges from the same cycle were deleted.
Thus, there is still a $u_{B}-v_{B}$ path in $G^{\prime}$ consisting of the other edges in the cycle.
But, since $u_{B} \in W$, $v_{B} \notin W$ and $e_{W}$ is the unique edge between $W$ and $G^{\prime} \setminus W$,
it must be that $e_{W}$ belongs to this path,
i.e.~$e_{W}$ must have been one of the other (at most $5$) edges in the cycle.
Thus, we have at most $5 \lceil \alpha k \rceil$ possibilities for $e_{W}$,
and hence for $W$
(since $e_{W}$ must be `oriented' so that $v_{B} \notin W$).

Thus, if case (b) was used,
then $f_{H}(G^{\prime}) \leq (h+7) \lceil \alpha k \rceil$,
and so we have at most
$\left( ^{(h+7) \lceil \alpha k \rceil}_{\phantom{ww} \lceil \alpha k \rceil} \right)
(6^{2}+6^{3}+6^{4}+6^{5})^{\lceil \alpha k \rceil}
\leq ((h+7)(6^{2}+ 6^{3}+ 6^{4}+ 6^{5})e)^{\lceil \alpha k \rceil}$
possibilities for $G$. \\

We have shown that each graph in $\mathcal{P}((1+\delta)k,d,D)$ is constructed at most
$((h+7)e)^{\lceil \alpha k \rceil} + ((h+7)(6^{2}+6^{3}+6^{4}+6^{5})e)^{\lceil \alpha k \rceil}
\leq 2 ((h+7)(6^{2}+6^{3}+6^{4}+6^{5})e)^{\lceil \alpha k \rceil}
\leq y^{\lceil \alpha k \rceil}$ times,
where $y$ denotes $2e(h+d+2)(6^{2}+6^{3}+6^{4}+6^{5})$.
Thus, the number of \textit{distinct} graphs
that we have constructed in $\mathcal{P}((1+\delta)k,d,D)$ must be at least
\begin{eqnarray*}
& & e^{-\alpha k} (1-\epsilon)^{k} ((1+\delta)k)! \left( \gamma_{d,D} \right)^{k}
\left( \frac{p}{172h! \alpha y} \right)^{\lceil \alpha k \rceil} \\
& \geq & e^{-\alpha k} (1-\epsilon)^{k} ((1+\delta)k)! \left( \gamma_{d,D} \right)^{(1+\delta)k}
\left( \gamma_{d,D} \right)^{- \lceil \alpha k \rceil h}
\left( \frac{p}{172h! \alpha y} \right)^{\lceil \alpha k \rceil}, \\
& & \textrm{ since } \delta k = \lceil \alpha k \rceil h \\
& \geq & (1-\epsilon)^{k} ((1+\delta)k)! \left( \gamma_{d,D} \right)^{(1+\delta)k}
\left( \frac{172 e h! \alpha y \left( \gamma_{d,D} \right)^{h}}{p} \right)^{- \lceil \alpha k \rceil}, \\
& & \textrm{ since } e^{-\alpha k} \geq e^{- \lceil \alpha k \rceil} \\
& \geq & (1-\epsilon)^{k} ((1+\delta)k)! \left( \gamma_{d,D} \right)^{(1+\delta)k}
(\alpha \beta)^{- \lceil \alpha k \rceil} \\
& \geq & (1-\epsilon)^{k} ((1+\delta)k)! \left( \gamma_{d,D} \right)^{(1+\delta)k}
(\alpha \beta)^{- \alpha k}
\textrm{, \phantom{www} since } \alpha \beta < 1 \\
& = & \left( \frac{1-\epsilon}{1-3\epsilon} \right)^{k} ((1+\delta)k)! \left( \gamma_{d,D} \right)^{(1+\delta)k}
\textrm{, \phantom{www} since } (\alpha \beta)^{\alpha} = 1-3 \epsilon \\
& \geq & \left( \frac{1-\epsilon}{1-3\epsilon} \right)^{k}
\frac{|\mathcal{P}((1+\delta)k,d,D)|}{(1+\epsilon)^{(1+\delta)k}}
\textrm{, \phantom{www} by (\ref{eq:gr})} \\
& > & \left( \frac{1-\epsilon}{(1-3\epsilon)(1+\epsilon)^{2}} \right)^{k}
|\mathcal{P}((1+\delta)k,d,D)|, \phantom{w} \textrm{ since } \delta < 1 \textrm{ (page~\pageref{delta})} \\
& > & |\mathcal{P}((1+\delta)k,d,D)|, \phantom{www} \textrm{ since }
(1-3\epsilon)(1+\epsilon)^{2} = 1 - \epsilon - 5 \epsilon^{2} - 3 \epsilon^{3}.
\end{eqnarray*}
Hence, we have our desired contradiction.
$\phantom{qwerty}$
\setlength{\unitlength}{0.25cm}
\begin{picture}(1,1)
\put(0,0){\line(1,0){1}}
\put(0,0){\line(0,1){1}}
\put(1,1){\line(-1,0){1}}
\put(1,1){\line(0,-1){1}}
\end{picture} \\
\\

As mentioned, we shall now see that we can actually drop the conditions that $d(n)$ is a constant
and $D(n)$ is monotonically non-decreasing:

\begin{Theorem}\label{bounded11}
Let $H$ be a (fixed) connected planar graph on $\{ 1,2, \ldots, h \}$.
Then there exists a constant $\alpha (h) >0$ such that,
given any integer-valued functions $d(n)$ and $D(n)$ satisfying
$\limsup_{n \to \infty}d(n) \leq \delta(H)$
and $\liminf_{n \to \infty}D(n) \geq \max \{\Delta(H), \deg_{H}(1)+1, 3 \},$
we have
\begin{displaymath}
\mathbf{P} [ f_{H} (P_{n,d,D}) \leq \alpha n ] < e^{-\alpha n}
\textrm{ for all sufficiently large } n.
\end{displaymath}
\end{Theorem}
\textbf{Proof}
Suppose we can find a graph $H$ and functions $d(n)$ and $D(n)$
that satisfy the conditions of this lemma,
but not the conclusion,
and let $\alpha = \alpha(h)$ be as given by Lemma~\ref{bounded106}.
Then there exist arbitrarily large `bad' $n$ for which
$\mathbf{P} [ f_{H} (P_{n,d,D}) \leq \alpha n ] \geq e^{-\alpha n}$.

Let $n_{1}$ be one of these bad $n$
and let us try to find a bad $n_{2}>n_{1}$ with $D(n_{2}) \geq D(n_{1})$.
Let us then try to find a bad $n_{3}>n_{2}$ with $D(n_{3}) \geq D(n_{2})$,
and so on.
We will either
(a) obtain an infinite sequence $n_{1},n_{2},n_{3} \ldots$
with $n_{1}<n_{2}<n_{3}< \ldots$
and $D(n_{1}) \leq D(n_{2}) \leq D(n_{3}) \leq \ldots$,
or (b) we will find a value $n_{k}$ such that all bad $n>n_{k}$ have $D(n) \leq D(n_{k})$.

Note that we must have $d(n) \in \{0,1,2,3,4,5\}$ for all $n$.
Hence, in case (a) there must exist a constant $d$ such that infinitely many of our $n_{i}$ satisfy $d(n_{i})=d$
(we shall call these $n_{i}$ `special').
Let the function $D^{*}$ be defined by setting $D^{*}(n)=D(n_{1})$ for all $n \leq n_{1}$
and $D^{*}(n) = D(n_{j})$ for all $n \in \{ n_{j-1}+1,n_{j-1}+2, \ldots, n_{j} \}$ for all $j> 1$.
Then $D^{*}$ is a monotonically non-decreasing integer-valued function satisfying
$\liminf_{n \to \infty} D^{*}(n) \geq \liminf_{n \to \infty} D(n)
\geq \max \{ \Delta(H), \deg_{H}(1)+1, 3 \}$.
Hence, since $d \leq \limsup_{n \to \infty} d(n) \leq \delta(H)$,
by Lemma~\ref{bounded106}
it must be that we have
$\mathbf{P} [ f_{H}(P_{n,d,D^{*}}) \leq \alpha n ] < e^{-\alpha n}$
for all sufficiently large $n$.
But recall that our infinitely many `special' $n_{i}$ satisfy
$(d, D^{*}(n_{i})) = (d(n_{i}), D(n_{i}))$,
and so
$\mathbf{P} [ f_{H}(P_{n_{i},d,D^{*}}) \leq \alpha n_{i} ]
= \mathbf{P} [ f_{H} (P_{n_{i},d,D}) \leq \alpha n_{i} ] \geq e^{-\alpha _{i}}$
for these $n_{i}$.
Thus, we obtain a contradiction.

In case (b),
note that we have $d(n_{i}) \in \{ 0,1,2,3,4,5 \}$ for all $i$
and that we also have
$D(n_{i}) \in \{ 3,4, \ldots, D(n_{k}) \}$ for all $i \geq k$.
Hence, there must exist constants $d$ and $D$ such that infintely many of our $n_{i}$ satisfy
$(d(n_{i}),D(n_{i})) = (d,D)$.
But by Lemma~\ref{bounded106} we have
$\mathbf{P} [ f_{H}(P_{n,d,D}) \leq \alpha n ] < e^{-\alpha n}$
for all large $n$,
and so we again obtain a contradiction.
$\phantom{qwerty}$
\setlength{\unitlength}{0.25cm}
\begin{picture}(1,1)
\put(0,0){\line(1,0){1}}
\put(0,0){\line(0,1){1}}
\put(1,1){\line(-1,0){1}}
\put(1,1){\line(0,-1){1}}
\end{picture} \\
\\

In the remainder of this section,
we will look at the case when $d(n) = D(n)$.
This time, we shall find it more convenient to introduce the concept of `$2$-appearances':

\begin{Definition} \label{2appdefn}
Let $J$ be a connected graph on the vertices $\{ 1,2, \ldots, |J| \}$.
Given a graph $G$,
we say that $J$ \emph{2-appears} at $W \subset V(G)$ if
(a) the increasing bijection from $\{ 1,2, \ldots, |J| \}$ to $W$ gives an isomorphism
between $J$ and the induced subgraph $G[W]$ of $G$;
and (b) there are exactly two edges, $e_{1} = r_{1}v_{1}$ and $e_{2} = r_{2}v_{2}$,
in $G$ between $W \supset \{ r_{1},r_{2} \}$ and $V(G) \setminus W \supset \{ v_{1}, v_{2} \}$,
these edges are non-adjacent
(i.e.~$r_{1} \neq r_{2}$ and $v_{1} \neq v_{2}$),
and $v_{1}$ and $v_{2}$ are also non-adjacent
(see Figure~\ref{2appfig}).

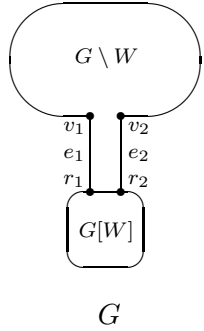
\begin{figure} [ht]
\setlength{\unitlength}{1cm}
\begin{picture}(20,4.25)(-5,-0.75)

\put(0.75,0){\line(1,0){0.5}}
\put(0.75,3.5){\line(1,0){0.5}}
\put(0.75,1){\line(1,0){0.5}}

\put(0.8,1){\line(0,1){1}}
\put(1.2,1){\line(0,1){1}}

\put(0.75,0.5){\oval(0.5,1)[l]}
\put(0.75,2.75){\oval(2,1.5)[l]}
\put(1.25,0.5){\oval(0.5,1)[r]}
\put(1.25,2.75){\oval(2,1.5)[r]}

\put(0.8,1){\circle*{0.1}}
\put(0.8,2){\circle*{0.1}}
\put(1.2,1){\circle*{0.1}}
\put(1.2,2){\circle*{0.1}}

\put(0.6,2.7){$G \setminus W$}
\put(0.65,0.4){\small{$G[W]$}}

\put(0.9,-0.75){\large{$G$}}

\put(1.3,1.45){$e_{2}$}
\put(1.3,1.8){$v_{2}$}
\put(1.3,1.1){$r_{2}$}

\put(0.45,1.45){$e_{1}$}
\put(0.45,1.8){$v_{1}$}
\put(0.45,1.1){$r_{1}$}

\end{picture}

\caption{$\textrm{A $2$-appearance at $W$ in $G$.}$} \label{2appfig}
\end{figure}
\end{Definition}

We will now give our result on $2$-appearances in $P_{n,d,D}$.
As the method is very similar to that of Lemma~\ref{bounded106},
we omit the full details (see Lemma~87 of~\cite{dow} for a complete proof):

\begin{Theorem} \label{bounded110}
Let $D \geq 3$ be a constant,
let $H$ be a (fixed) $D$-regular connected planar graph on $\{1,2, \ldots,h \}$,
and let $f \in E(H)$ be a non cut-edge.
Then there exist constants $\alpha(h) > 0$
and $N(h)$ such that
\begin{eqnarray*}
\mathbf{P}[P_{n,D,D}
\textrm{ will not have a set of at least $\alpha n$ $2$-appearances of $H \setminus f$}] \\
< e^{- \alpha n}
\left \{ \begin{array}{ll}
\textrm{for all } n \geq N \textrm{ if } D=4 \\
\textrm{for all even } n \geq N \textrm{ if } D \in \{ 3,5 \}.
\end{array} \right.
\end{eqnarray*}
\end{Theorem}
\textbf{Sketch of Proof}
We choose a specific $\alpha$
and suppose that the result is false for $n=k$,
where $k$ is suitably large.
Using Theorem~\ref{bounded104}/Theorem~\ref{bounded105}
(depending on the parity of $D$),
it then follows that there are many graphs $G \in \mathcal{P}(k,D,D)$ with at most $\alpha k$
$2$-appearances of $H \setminus f$.

From each such $G$,
we construct graphs in $\mathcal{P}((1+\delta)k,D,D)$, for a fixed $\delta >0$,
by replacing some edges in $G$ with $2$-appearances of $H \setminus f$.

The fact that the original graphs had few $2$-appearances of $H \setminus f$
can then be used to show that there is not much double-counting,
and so we find that we have built so many graphs in $\mathcal{P}((1+\delta)k,D,D)$
that we contradict Theorem~\ref{bounded104}/Theorem~\ref{bounded105}. \\
\\

\section{Components} \label{cpts}

We shall now use our appearance results from the previous section
to investigate the probability of $P_{n,d,D}$ having given components.

We already know from Section~\ref{bconn} that
(assuming $D(n) \geq 3$ for all $n$, as always)
$\liminf \mathbf{P}[P_{n,d,D}$ will be connected$] >0$,
so certainly it must be that
$\limsup \mathbf{P}[P_{n,d,D}$ will have a component isomorphic to $H] <1$ for all $H$.
In this section,
we will now see (in Theorem~\ref{bounded404}) that for all feasible $H$ we also have
$\liminf \mathbf{P}[P_{n,d,D}$ will have a component isomorphic to $H] >0$.

As we are going to be using Theorems~\ref{bounded11} and~\ref{bounded110} from Section~\ref{apps},
we will start by dealing with the $d(n) < D(n)$ and $d(n)=D(n)$ cases separately
(in Lemmas~\ref{bounded7} and~\ref{bounded111}, respectively),
but we shall then combine these results in Theorem~\ref{bounded404}. \\

We start with the case when $d(n) < D(n)$ for all $n$:

\begin{Lemma} \label{bounded7}
Let $d(n)$ and $D(n)$ be any integer-valued functions
that for all $n$ satisfy $D(n) \geq 3$ and $d(n)<D(n)$.
Then, given any (fixed) connected planar graph $H$ with
$\limsup_{n \to \infty}d(n) \leq \delta(H) \leq \Delta(H) \leq \liminf_{n \to \infty}D(n)$,
we have
\begin{displaymath}
\liminf_{n \to \infty}
\mathbf{P}[P_{n,d,D}
\textrm{ will have a component isomorphic to $H$}] > 0.
\end{displaymath}
\end{Lemma}
\textbf{Proof}
Let $\mathcal{N}^{<}$ denote the set of values of $n$ for which $\delta(H)<D(n)$.
We shall start by proving the result for $\mathcal{N}^{<}$.

Without loss of generality (by symmetry),
we may assume that we have
$V(H) = \{ 1,2, \ldots, |H| \}$
and that $\deg_{H}(1) = \delta (H)$.
Thus, by Theorem \ref{bounded11}, we know
there exist $\alpha >0$ and $N_{1}$ such that
$\mathbf{P}[f_{H}(P_{n,d,D}) \leq \alpha n] < e^{- \alpha n}$
for $\{ n \in \mathcal{N}^{<} : n \geq N_{1} \}$.
Hence, given any $\epsilon > 0$,
there certainly exists an $N_{2}$
such that $\mathbf{P}[f_{H}(P_{n,d,D}) \geq \alpha n] > \epsilon$
for $\{ n \in \mathcal{N}^{<} : n \geq N_{2} \}$.
Let us suppose that we can find an $n \in \mathcal{N}^{<}$ satisfying $n \geq N_{2}$ for which
$\mathbf{P}[P_{n,d,D}$
will have a component isomorphic to $H] < \frac{\epsilon}{2}$
(if not, then we are done)
and let $\mathcal{G}(n,d,D)$ denote the set of graphs in $\mathcal{P}(n,d,D)$
which have both
(i) at least $\alpha n$ cut-able appearances of $H$
and (ii) no components isomorphic to $H$.

\begin{figure} [ht]
\setlength{\unitlength}{1cm}
\begin{picture}(20,3.5)(-2.25,0)

\put(0.75,0){\line(1,0){0.5}}
\put(0.75,3.5){\line(1,0){0.5}}

\put(0.75,1){\line(1,0){0.5}}
\put(0.75,2){\line(1,0){0.5}}

\put(0.75,0.5){\oval(0.5,1)[l]}
\put(0.75,2.75){\oval(2,1.5)[l]}
\put(1.25,0.5){\oval(0.5,1)[r]}
\put(1.25,2.75){\oval(2,1.5)[r]}

\put(1,1){\circle*{0.1}}
\put(1,2){\circle*{0.1}}

\put(1,1){\line(0,1){1}}

\put(3.25,1.75){\vector(1,0){1}}

\put(6.25,0){\line(1,0){0.5}}
\put(6.25,3.5){\line(1,0){0.5}}

\put(6.25,1){\line(1,0){0.5}}
\put(6.25,2){\line(1,0){0.5}}

\put(6.25,0.5){\oval(0.5,1)[l]}
\put(6.25,2.75){\oval(2,1.5)[l]}
\put(6.75,0.5){\oval(0.5,1)[r]}
\put(6.75,2.75){\oval(2,1.5)[r]}

\put(6.5,1){\circle*{0.1}}
\put(6.5,2){\circle*{0.1}}

\end{picture}

\caption{Using an appearance to construct a component isomorphic to $H$.} \label{bcptfig2}
\end{figure}
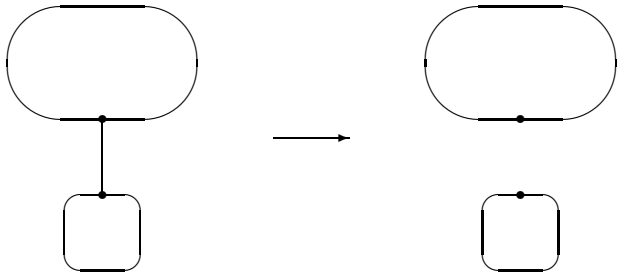

We may construct a graph in $\mathcal{P}(n,d,D)$
with a component isomorphic to $H$ simply by taking a graph in $\mathcal{G}(n,d,D)$
(at least $\frac{\epsilon}{2} |\mathcal{P}(n,d,D)|$ choices)
and deleting the associated cut-edge from a cut-able appearance of $H$
(at least $\alpha n$ choices).
Thus,
we have at least
$\frac{\epsilon \alpha n}{2} |\mathcal{P}(n,d,D)|$
ways to construct (not necessarily distinct) graphs in $\mathcal{P}(n,d,D)$
with a component isomorphic to $H$
(see Figure~\ref{bcptfig2}).

Given one of our constructed graphs,
there will be at most two possibilities for which is our deliberately created component
(since there were no components isomorphic to $H$ in the original graph
and we can have only created at most two when we deleted the cut-edge)
and then at most $n$ possibilities for the vertex in the rest of the graph
that this component was attached to originally.
Hence,
the amount of double-counting is at most $2n$,
and so we find that we can construct at least
$\frac{\epsilon \alpha}{4} |\mathcal{P}(n,d,D)|$
\textit{distinct} graphs in $\mathcal{P}(n,d,D)$
with a component isomorphic to $H$,
which is what we wanted to prove.

Now let $\mathcal{N}^{=}$ denote the set of values of $n$ for which $\delta(H)=D(n)$.
It only remains to prove the result for $\mathcal{N}^{=}$.
If $\mathcal{N}^{=} \neq \emptyset$,
then we must have $\delta(H) \geq 3$,
since $D(n) \geq 3$.
Hence, there exists a non cut-edge $f \in E(H)$.
Let $H^{\prime} = H \setminus f$
and note that $\delta(H^{\prime}) = D(n)-1 \geq d(n)$ for all $n \in \mathcal{N}^{=}$.
Thus, we may use Theorem~\ref{bounded11} with $H^{\prime}$
and then follow the same proof as with $\mathcal{N}^{<}$ to find constants
$\epsilon^{\prime} > 0$ and $N_{3}$ such that we have
$\mathbf{P}[P_{n,d,D}$
will have a component isomorphic to $H^{\prime}] > \epsilon^{\prime}$
for $\{ n \in \mathcal{N}^{=} : n \geq N_{3} \}$.
Let us suppose that we can find an $n \in \mathcal{N}^{=}$ satisfying $n \geq N_{3}$ for which
$\mathbf{P}[P_{n,d,D}
\textrm{ will have a component isomorphic to }H] < \frac{\epsilon^{\prime}}{2}$
(if not, then we are done)
and let $\mathcal{H}(n,d,D)$ denote the set of graphs in
$\mathcal{P}(n,d,D)$
with a component isomorphic to $H^{\prime}$ but without any components isomorphic to $H$.

Given a graph in $\mathcal{H}(n,d,D)$,
we may construct a graph in $\mathcal{P}(n,d,D)$
with a component isomorphic to $H$ simply by choosing a component isomorphic to $H^{\prime}$
and adding an appropriate edge.
The amount of double-counting will be at most $3|H|$,
since we will know exactly where the modified component is,
and so we find that the number of \textit{distinct} graphs in $\mathcal{P}(n,d,D)$
with a component isomorphic to $H$ must be at least
$\frac{|\mathcal{H}(n,d,D)|}{3|H|} \geq \frac{\epsilon^{\prime}|\mathcal{P}(n,d,D)|}{6|H|}$,
and so we are done.~
\phantom{qwerty}
$\setlength{\unitlength}{.25cm}
\begin{picture}(1,1)
\put(0,0){\line(1,0){1}}
\put(0,0){\line(0,1){1}}
\put(1,1){\line(-1,0){1}}
\put(1,1){\line(0,-1){1}}
\end{picture}$ \\
\\

We shall now see an analogous result for when $d(n)=D(n)$ for all $n$:

\begin{Lemma} \label{bounded111}
Let $D \geq 3$ be a constant
and let $H$ be a (fixed) $D\textrm{-regular}$ connected planar graph.
Then there exist constants $\epsilon(H) > 0$ and $N(H)$ such that
\begin{eqnarray*}
\mathbf{P}[P_{n,D,D}
\textrm{ will have a component isomorphic to } H] \\
> \epsilon
\left \{ \begin{array}{ll}
\textrm{for all } n \geq N \textrm{ if } D=4 \\
\textrm{for all even } n \geq N \textrm{ if } D \in \{ 3,5 \}.
\end{array} \right.
\end{eqnarray*}
\end{Lemma}
\textbf{Proof}
In order to simplify parity matters,
we shall just prove the result for $D=4$,
but the $D \in \{ 3,5 \}$ cases will follow in a completely analogous way.

Without loss of generality,
we may assume that $V(H) = \{ 1,2, \ldots, |H| \}$.
Let $f \in E(H)$ be an arbitrary non cut-edge.
Then, by Theorem \ref{bounded110},
we know that there exists $\alpha >0$ and there exists $N_{1}$ such that for all $n \geq N_{1}$ we have
$\mathbf{P}[P_{n,D,D}$ does not have
$\geq \alpha n$ $2$-appearances of $H \setminus f] < e^{- \alpha n}$.
Thus, given any $\delta > 0$,
there certainly exists an $N_{2}$
such that $\mathbf{P}[P_{n,d,D}$ has
$\geq \alpha n$ $2$-appearances of $H \setminus f] > \delta$ for all $n \geq N_{2}$.
Let us now suppose that we can find an $n \geq N_{2}$ for which
$\mathbf{P}[P_{n,d,D}
\textrm{ will have a component isomorphic to }H] < \frac{\delta}{2}$
(if not, then we are done)
and let $\mathcal{G}(n,d,D)$ denote the set of graphs in $\mathcal{P}(n,d,D)$
which have both
(i) at least $\alpha n$ $2$-appearances of $H$
and (ii) no components isomorphic to $H$.

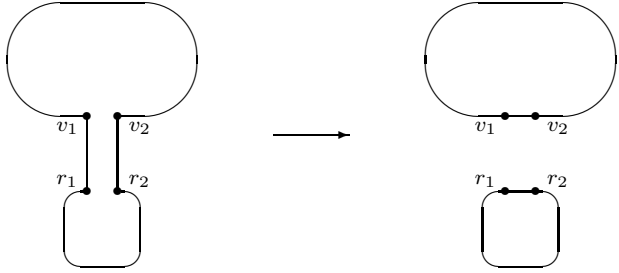
\begin{figure} [ht]
\setlength{\unitlength}{1cm}
\begin{picture}(20,3.5)(-2.25,0)

\put(0.75,0){\line(1,0){0.5}}
\put(0.75,3.5){\line(1,0){0.5}}

\put(0.8,1){\line(0,1){1}}
\put(1.2,1){\line(0,1){1}}

\put(0.75,0.5){\oval(0.5,1)[l]}
\put(0.75,2.75){\oval(2,1.5)[l]}
\put(1.25,0.5){\oval(0.5,1)[r]}
\put(1.25,2.75){\oval(2,1.5)[r]}

\put(0.8,1){\circle*{0.1}}
\put(0.8,2){\circle*{0.1}}
\put(1.2,1){\circle*{0.1}}
\put(1.2,2){\circle*{0.1}}

\put(1.35,1.8){$v_{2}$}
\put(1.35,1.1){$r_{2}$}

\put(0.4,1.8){$v_{1}$}
\put(0.4,1.1){$r_{1}$}

\put(3.25,1.75){\vector(1,0){1}}

\put(6.25,0){\line(1,0){0.5}}
\put(6.25,3.5){\line(1,0){0.5}}

\put(6.25,1){\line(1,0){0.5}}
\put(6.25,2){\line(1,0){0.5}}

\put(6.25,0.5){\oval(0.5,1)[l]}
\put(6.25,2.75){\oval(2,1.5)[l]}
\put(6.75,0.5){\oval(0.5,1)[r]}
\put(6.75,2.75){\oval(2,1.5)[r]}

\put(6.3,1){\circle*{0.1}}
\put(6.3,2){\circle*{0.1}}
\put(6.7,1){\circle*{0.1}}
\put(6.7,2){\circle*{0.1}}

\put(6.85,1.8){$v_{2}$}
\put(6.85,1.1){$r_{2}$}

\put(5.9,1.8){$v_{1}$}
\put(5.9,1.1){$r_{1}$}

\end{picture}

\caption{Using a $2$-appearances to construct a component isomorphic to $H$.} \label{bcptfig}
\end{figure}

We may construct a graph in $\mathcal{P}(n,d,D)$
with a component isomorphic to $H$ simply by taking a graph $G \in \mathcal{G}(n,d,D)$
(at least $\frac{\delta}{2} |\mathcal{P}(n,d,D)|$ choices);
choosing a $2$-appearance $W$ of $H \setminus e$
(at least $\beta n$ choices);
deleting the two edges of the form $v_{1}r_{1}$ and $v_{2}r_{2}$ for $\{ r_{1}, r_{2} \} \subset W$
and $\{ v_{1}, v_{2} \} \subset V(G) \setminus W$;
and then finally inserting the two edges $v_{1}v_{2}$ and $r_{1}r_{2}$
(note that $v_{1}$ and $v_{2}$ were not originally adjacent,
by the definition of a $2$-appearance,
and that $r_{1}$ and $r_{2}$ were also not originally adjacent,
since they must be the two vertices of degree $D-1$ in $G[W]$,
by $D$-regularity of $G$,
and we know $G[W]$ is isomorphic to $H \setminus f$).
Thus,
we have at least $\frac{\delta \beta n}{2} |\mathcal{P}(n,d,D)|$
ways to construct (not necessarily distinct) graphs in $\mathcal{P}(n,d,D)$
with a component isomorphic to $H$
(see Figure~\ref{bcptfig}).

Given one of our constructed graphs,
there will be at most two possibilities for which is our deliberately constructed component
(since there were no components isomorphic to $H$ in the original graph
and we can have only created at most two).
We then know which edge was inserted into it and have at most $3n$ (by planarity) possibilities for
which edge was inserted in the rest of the graph.
We also have two further possibilities for how the vertices in these two edges were connected originally.
Hence,
the amount of double-counting is at most $12n$,
and so we find that we can construct at least
$\frac{\delta \beta}{12} |\mathcal{P}(n,d,D)|$
\textit{distinct} graphs in $\mathcal{P}(n,d,D)$
with a component isomorphic to $H$,
which is what we wanted to prove.~
\phantom{qwerty}
\begin{picture}(1,1)
\put(0,0){\line(1,0){1}}
\put(0,0){\line(0,1){1}}
\put(1,1){\line(-1,0){1}}
\put(1,1){\line(0,-1){1}}
\end{picture} \\
\\

We may now combine Lemmas~\ref{bounded7} and~\ref{bounded111} to obtain our full result:

\begin{Theorem} \label{bounded404}
Let $d(n)$ and $D(n)$ be any integer-valued functions, subject to $D(n) \geq 3$ for all $n$
and $(d(n),D(n)) \notin \{ (3,3),(5,5) \}$ for odd $n$.
Then, given any (fixed) connected planar graph
$H$ satisfying
$\limsup_{n \to \infty} d(n) \leq \delta(H) \leq \Delta(H)
\leq \liminf_{n \to \infty} D(n)$,
we have
\begin{displaymath}
\liminf_{n \to \infty}\mathbf{P}[P_{n,d,D}
\textrm{ will have a component isomorphic to $H$}] > 0.
\end{displaymath} \\
\end{Theorem}

\section{Subgraphs} \label{subs}

We will now use the results of the previous two sections to investigate the probability of $P_{n,d,D}$
having \textit{subgraphs} isomorphic to given connected planar graphs.
As always, we shall assume throughout that $D(n) \geq 3$ for all $n$.

Clearly, for those values of $n$ for which $D(n) < \Delta(H)$,
we must have $\mathbf{P}[P_{n,d,D} \textrm{ will have a copy of } H] =0$.
For sufficiently large $n$,
it turns out that the only other time when we can have this
is if $d(n)=D(n)=4$
and $H$ happens to be a graph that can \textit{never} be a subgraph of a $4$-regular planar graph.

Apart from the above exceptions,
we shall see that the matter of whether
$\mathbf{P}[P_{n,d,D} \textrm{ will have a copy of } H]$
is bounded away from $0$ and/or $1$
actually depends only on whether or not $H$ is $D(n)$-regular.
For those values of $n$ for which this is the case,
any copy of $H$ must be a component and so it already follows
from Sections~\ref{bconn} and~\ref{cpts} that
the probability
must indeed be bounded away from both $0$ and $1$.
If there aren't arbitrarily large values of $n$ for which $H$ is $D(n)$-regular,
though,
we shall be able to use our appearance results of Section~\ref{apps}
to deduce (in Theorem~\ref{bounded1001}) that
$\mathbf{P}[P_{n,d,D} \textrm{ will have a copy of } H] \to 1$
(again, with the exception of the cases given in the previous paragraph). \\

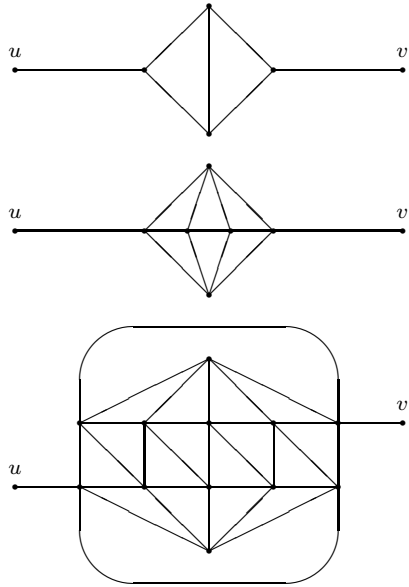
\begin{figure} [ht]
\setlength{\unitlength}{0.85cm}
\begin{picture}(20,9)(-6.1,-6)

\put(-2,2){\line(1,0){2}}
\put(2,2){\line(1,0){2}}

\put(1,3){\line(-1,-1){1}}
\put(1,3){\line(1,-1){1}}

\put(1,1){\line(-1,1){1}}
\put(1,1){\line(1,1){1}}

\put(1,3){\line(0,-2){2}}

\put(0,2){\circle*{0.1}}
\put(2,2){\circle*{0.1}}
\put(1,3){\circle*{0.1}}
\put(1,1){\circle*{0.1}}

\put(-2,2){\circle*{0.1}}
\put(4,2){\circle*{0.1}}

\put(-2.1,2.2){$u$}
\put(3.9,2.2){$v$}

\put(-2,-0.5){\line(1,0){6}}

\put(1,0.5){\line(-1,-1){1}}
\put(1,0.5){\line(-1,-3){0.333333}}
\put(1,0.5){\line(1,-3){0.333333}}
\put(1,0.5){\line(1,-1){1}}

\put(1,-1.5){\line(-1,1){1}}
\put(1,-1.5){\line(-1,3){0.333333}}
\put(1,-1.5){\line(1,3){0.333333}}
\put(1,-1.5){\line(1,1){1}}

\put(0,-0.5){\circle*{0.1}}
\put(0.666667,-0.5){\circle*{0.1}}
\put(1.333333,-0.5){\circle*{0.1}}
\put(2,-0.5){\circle*{0.1}}
\put(1,0.5){\circle*{0.1}}
\put(1,-1.5){\circle*{0.1}}

\put(-2,-0.5){\circle*{0.1}}
\put(4,-0.5){\circle*{0.1}}

\put(-2.1,-0.3){$u$}
\put(3.9,-0.3){$v$}

\put(-1,-3.5){\line(1,0){4}}
\put(-1,-4.5){\line(1,0){4}}

\put(-1,-3.5){\line(0,-1){1}}
\put(0,-3.5){\line(0,-1){1}}
\put(1,-2.5){\line(0,-1){3}}
\put(2,-3.5){\line(0,-1){1}}
\put(3,-3.5){\line(0,-1){1}}

\put(1,-2.5){\line(-2,-1){2}}
\put(1,-2.5){\line(-1,-1){1}}
\put(1,-2.5){\line(1,-1){1}}
\put(1,-2.5){\line(2,-1){2}}

\put(1,-5.5){\line(-2,1){2}}
\put(1,-5.5){\line(-1,1){1}}
\put(1,-5.5){\line(1,1){1}}
\put(1,-5.5){\line(2,1){2}}

\put(1,-4.5){\oval(4,3)[b]}
\put(1,-3.5){\oval(4,3)[t]}

\put(-1,-3.5){\line(1,-1){1}}
\put(0,-3.5){\line(1,-1){1}}
\put(1,-3.5){\line(1,-1){1}}
\put(2,-3.5){\line(1,-1){1}}
\put(3,-3.5){\line(1,0){1}}

\put(-1,-3.5){\circle*{0.1}}
\put(0,-3.5){\circle*{0.1}}
\put(1,-3.5){\circle*{0.1}}
\put(2,-3.5){\circle*{0.1}}
\put(3,-3.5){\circle*{0.1}}

\put(-1,-4.5){\circle*{0.1}}
\put(0,-4.5){\circle*{0.1}}
\put(1,-4.5){\circle*{0.1}}
\put(2,-4.5){\circle*{0.1}}
\put(3,-4.5){\circle*{0.1}}

\put(1,-2.5){\circle*{0.1}}
\put(1,-5.5){\circle*{0.1}}

\put(-2,-4.5){\line(1,0){1}}

\put(-2.1,-4.3){$u$}
\put(3.9,-3.3){$v$}

\put(4,-3.5){\circle*{0.1}}
\put(-2,-4.5){\circle*{0.1}}

\end{picture}

\caption{Replacing the edge $f=uv$ (cases $D=3$, $D=4$ and $D=5$).} \label{edgefig}
\end{figure}

We will start by working towards the $d(n)=D(n)$ case:

\begin{Lemma} \label{bounded901}
Let $H$ be a (fixed) connected planar graph
and let $D \in \{3,4,5 \}$ be a constant.
Suppose $H$ is not $D$-regular,
but that there exists a $D$-regular planar graph $H^{*}$ that contains a copy of $H$.
Then there exist constants $\beta(H) > 0$ and $N(H)$ such that
\begin{eqnarray*}
\mathbf{P}[P_{n,D,D}
\textrm{ will \emph{not} have at least $\beta n$ copies of $H$}]
< e^{- \beta n} \\
\left \{ \begin{array}{ll}
\textrm{for all } n \geq N \textrm{ if } D=4 \\
\textrm{for all even } n \geq N \textrm{ if } D \in \{ 3,5 \}.
\end{array} \right.
\end{eqnarray*}
\end{Lemma}
\textbf{Proof}
Since $H$ is not $D$-regular,
it must be that $H^{*}$ contains an edge $f=uv$
such that $H^{*} \setminus f$ also contains a copy of $H$.
Without loss of generality,
we may assume that $f$ is not a cut-edge in $H^{*}$,
since we could replace $f_{i}$ with a copy of the appropriate graph from Figure~\ref{edgefig}
and in this way obtain a $D$-regular planar graph containing several non cut-edges
that don't interfere with our copy of $H_{i}$.
Thus, the result then follows from Theorem~\ref{bounded110} with $H^{*} \setminus f$.
$\phantom{qwerty}
\setlength{\unitlength}{.25cm}
\begin{picture}(1,1)
\put(0,0){\line(1,0){1}}
\put(0,0){\line(0,1){1}}
\put(1,1){\line(-1,0){1}}
\put(1,1){\line(0,-1){1}}
\end{picture}$ \\
\\

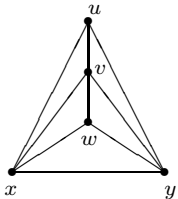
\begin{figure} [ht]
\setlength{\unitlength}{1cm}
\begin{picture}(10,2)(-5,-0.2)

\put(0,0){\line(1,0){2}}
\put(1,0.666667){\line(0,1){1.333333}}
\put(0,0){\line(3,2){1}}
\put(0,0){\line(3,4){1}}
\put(0,0){\line(1,2){1}}
\put(2,0){\line(-3,2){1}}
\put(2,0){\line(-3,4){1}}
\put(2,0){\line(-1,2){1}}

\put(0,0){\circle*{0.1}}
\put(2,0){\circle*{0.1}}
\put(1,0.666667){\circle*{0.1}}
\put(1,1.333333){\circle*{0.1}}
\put(1,2){\circle*{0.1}}

\put(-0.1,-0.3){$x$}
\put(2,-0.3){$y$}
\put(0.9,0.366667){$w$}
\put(1.08,1.3){$v$}
\put(1,2.1){$u$}

\end{picture}
\caption{The unique planar embedding of $K_{5} \setminus \{ u,w \}$.} \label{PH}
\end{figure}

If $D \in \{ 3,5 \}$,
then it can easily be shown that there does always exist a $D$-regular planar graph $H^{*} \supset H$
whenever $\Delta (H) \leq D$,
since in these cases there exist planar graphs
that are $D$-regular except for exactly one vertex with degree $D-1$,
and hence we can extend $H$ into a $D$-regular planar graph
simply by attaching an appropriate number of these graphs to any vertices of $H$ that have degree less than $D$.
This trick does not work for $D=4$, however,
since clearly a graph that is $4$-regular except for exactly one vertex of degree $3$ would have to have
an odd sum of degrees!
In fact,
the following example shows that
there do actually exist some planar graphs with maximum degree at most $4$
that can't ever be contained in a $4$-regular planar graph:

\begin{Example} \label{bounded908}
No $4\textrm{-regular}$ planar graph contains a copy of the graph $K_{5}$ minus an edge.
\end{Example}
\textbf{Proof}
The graph $K_{5} \setminus \{ u,w \}$ is drawn with its \textit{unique} planar embedding (see~\cite{whi})
in Figure~\ref{PH}.

Consider any planar graph $G \supset K_{5} \setminus \{ u,w \}$ with $\Delta (G) = 4$.
Since we already have
$\textrm{deg}_{H}(v) = \textrm{deg}_{H}(x) = \textrm{deg}_{H}(y) = 4$,
any new edge with at least one endpoint inside the triangle given by $vxy$
must have both endpoints inside.
Hence, the sum of degrees inside this triangle must remain odd,
and so this region must still contain a vertex of odd degree.
Thus, $G$ is not $4$-regular.
$\phantom{}
\setlength{\unitlength}{.25cm}
\begin{picture}(1,1)
\put(0,0){\line(1,0){1}}
\put(0,0){\line(0,1){1}}
\put(1,1){\line(-1,0){1}}
\put(1,1){\line(0,-1){1}}
\end{picture}$ \\

An $O\left( |H|^{2.5} \right)$ time algorithm for determining
whether or not a given graph $H$ can ever be a subgraph of a $4$-regular planar graph
is given in~\cite{add}. \\
\\

It now only remains for us to deal with the case when $d(n) < D(n)$.
But note that this can be deduced easily simply by applying Theorem~\ref{bounded11}
to a connected planar graph $H^{\prime} \supset H$
with sufficiently high minimum degree.
Hence, we have:

\begin{Lemma} \label{bounded904}
Let $H$ be a (fixed) connected planar graph.
Suppose $d(n)$ and $D(n)$ are integer-valued functions that for all large $n$ satisfy
(a) $d(n) < \min \{ 6, D(n) \}$, and
(b) $D(n) \geq \max \{ \Delta(H), \delta(H)+1, 3 \}$.
Then there exists a constant $\beta(H) > 0$ such that
\begin{displaymath}
\mathbf{P}[P_{n,d,D}
\textrm{ will \emph{not} have at least $\beta n$ copies of $H$}]
= e^{- \Omega (n)}.
\end{displaymath}
\end{Lemma}

\phantom{p}

Combining all the results of this section,
we obtain our full result:

\begin{Theorem} \label{bounded1001}
Let $H$ be a (fixed) connected planar graph.
Suppose $d(n)$ and $D(n)$ are integer-valued functions that for all $n$ satisfy
(a) $d(n) \leq \min \{ 5, D(n) \}$,
(b) $D(n) \geq \max \{ \Delta(H), \delta(H)+1, 3 \}$,
(c) $(d(n),D(n)) \notin \{ (3,3),(5,5) \}$ for odd $n$,
and also (d) $(d(n),D(n)) \neq (4,4)$
if $H$ happens to be a graph that can never be contained within a $4$-regular planar graph.
Then there exists a constant $\beta(H) > 0$ such that
\begin{displaymath}
\mathbf{P}[P_{n,d,D}
\textrm{ will \emph{not} have at least $\beta n$ copies of $H$}]
= e^{- \Omega (n)}.
\end{displaymath} \\
\end{Theorem}

\section*{Acknowledgements}

I am very grateful to Colin McDiarmid for his advice and helpful comments on the paper,
and also to an anonymous reviewer. \\
\\

\end{document}